\documentclass[12pt]{article}
\usepackage{amsmath,amsthm,amssymb,times,titlesec}
\usepackage{graphicx}
\usepackage{pdfsync}
\usepackage{float}
\numberwithin{equation}{section}
\textheight % 528 + 26 + 11 + 24 + 6 + 55 for luck
            650pt
\textwidth % 360 + 53 + 47 for luck
           460pt
\titleformat{\subsubsection}
{\normalfont\fontsize{13}{14}\selectfont\itshape}{\thesubsubsection}{1em}{}

\newcommand\scaleddot{\scalebox{.89}{.}}
\makeatletter
\renewcommand{\ddddot}[1]{%
  {\mathop{\kern\z@#1}\limits^{\makebox[0pt][c]{\vbox to-2\ex@{\kern-\tw@\ex@\hbox{\normalfont\scaleddot\kern-0.5pt\scaleddot\kern-0.5pt\scaleddot\kern-0.5pt\scaleddot}\vss}}}}}

\newcommand{\dr}{\, \mathrm{d}r}
\newcommand{\ds}{\, \mathrm{d}s}
\newcommand{\dt}{\, \mathrm{d}t}

\newcommand{\dz}{\, \mathrm{d}z}

\newcommand{\DD}{{\mathcal D}}

\newcommand{\LL}{{\mathcal L}}

\newcommand{\NN}{{\mathcal N}}

\newcommand{\PP}{{\mathcal P}}

\newcommand{\XX}{{\mathcal X}}

\DeclareMathOperator{\cosech}{cosech}

\newcommand{\sdfrac}[2]{\mbox{\small$\displaystyle\frac{#1}{#2}$}}

\newtheorem{theorem}{Theorem}[section]
\newtheorem{lemma}[theorem]{Lemma}

\newtheorem{definition}[theorem]{Definition}
\theoremstyle{definition}

\renewcommand{\i}{\mathrm{i}}

\newcounter{count}

\title{Spatial dynamics methods for solitary waves\\ on a ferrofluid jet}

\author{M. D. Groves\thanks{Fachrichtung 6.1 - Mathematik, Universit\"{a}t des Saarlandes,
Postfach 151150, 66041 Saarbr\"{u}cken, Germany; 
Department of Mathematical Sciences, Loughborough
University, Loughborough, LE11 3TU, UK
} \and
D. V. Nilsson\footnote{Centre for Mathematical Sciences, Lund University, PO Box 118, 22100 Lund, Sweden}}
\date{}

\begin{document}

\maketitle

\begin{abstract}

This paper presents existence theories for  several families of
axisymmetric solitary waves on the surface of an otherwise cylindrical ferrofluid
jet surrounding a stationary metal rod. The ferrofluid, which is governed by a general (nonlinear) magnetisation law, is
subject to an azimuthal magnetic field generated by an electric current flowing along the rod.

The ferrohydrodynamic problem for axisymmetric travelling waves is formulated as an infinite-dimensional Hamiltonian
system in which the axial direction is the
time-like variable. A centre-manifold reduction technique is employed
to reduce the system to a locally equivalent Hamiltonian system with a finite
number of degrees of freedom, and homoclinic solutions to the reduced system, which correspond
to solitary waves, are detected by dynamical-systems methods.
\end{abstract}

\section{Introduction}\label{Introduction}

\begin{figure}[h]
\centering
\includegraphics[scale=1.3]{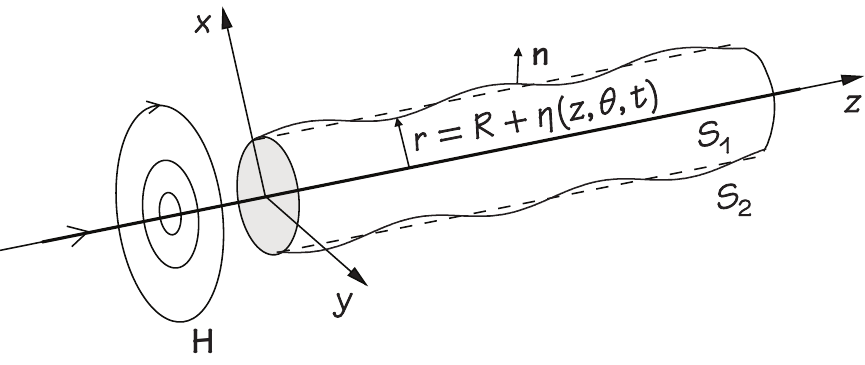}
{\it \caption{Waves on the surface of a ferrofluid jet surrounding a current-carrying wire.}
\label{Definition sketch}}
\end{figure}\pagebreak
 
We consider an incompressible, inviscid ferrofluid of unit density in the region
$$S_1:=\{0<r<R+\eta(\theta,z,t)\}$$
bounded by the free interface $\{r=R+\eta(\theta,z,t)\}$
and a current-carrying wire at $\{r=0\}$, where $(r,\theta,z)$ are cylindrical polar coordinates.
The fluid is subject to a static magnetic field and
the surrounding region $$S_2=\{r>R+\eta(\theta,z,t)\}$$
is a vacuum (see Figure \ref{Definition sketch}). \emph{Travelling waves}
move in the axial direction with constant speed $c$ and without
change of shape, so that $\eta(\theta,z,t) = \eta(\theta,z-ct)$. We are interested in particular
in \emph{axisymmetric solitary waves} for which $\eta$ does not depend upon $\theta$
and $\eta(z-ct) \rightarrow 0$ as $z-ct \rightarrow \pm \infty$.
Waves of this kind for ferrofluids with a linear magnetisation law have been
investigated using a weakly nonlinear approximation by Rannacher \& Engel \cite{RannacherEngel06},
experimentally by Bourdin, Bacri \& Falcon \cite{BourdinBacriFalcon10} and numerically by
Blyth \& Parau \cite{BlythParau14}. In this paper we present a rigorous existence theory
for small-amplitude solitary waves and consider fluids with a general (nonlinear) magnetisation law.

Our starting point is a formulation of the hydrodynamic problem as a reversible Hamiltonian system
\begin{equation}
\eta_z = \frac{\delta H}{\delta \omega}, \quad
\omega_z = -\frac{\delta H}{\delta \eta}, \quad
\hat{\phi}_z = \frac{\delta H}{\delta \hat{\zeta}}, \quad
\hat{\zeta}_z = -\frac{\delta H}{\delta \hat{\phi}}
\label{SHS}
\end{equation}
in which the axial coordinate $z$ plays the role of time, $\hat{\phi}$ is a variable related to the fluid velocity potential $\phi$ and
$\omega$, $\hat{\zeta}$ are the momenta associated with the coordinates $\eta$, $\hat{\phi}$. The spatial Hamiltonian
system \eqref{SHS} is derived from a variational principle for the governing equations in Section \ref{Spatial dynamics};
it depends upon two dimensionless physical parameters $\alpha$ and $\beta$ (see equation
\eqref{Defns of alpha and beta} for precise definitions) and the (dimensionless) magnitude $m_1(|{\bf H}_1|)$ of the magnetic intensity corresponding to the
magnetic field ${\bf H}_1$ in the ferrofluid.

\emph{Homoclinic solutions} of \eqref{SHS} (solutions with $(\eta,\omega,\hat{\phi},\hat{\zeta}) \rightarrow 0$ as $z \rightarrow \pm \infty$) are of particular interest since they correspond to solitary waves.
We detect such solutions using a technique known as the \emph{Kirchg\"{a}ssner reduction}
(Section \ref{Reduction}), in which
a centre-manifold reduction principle is used to show that all small, globally bounded solutions of
a spatial (Hamiltonian) evolutionary system solve a (Hamiltonian) system of ordinary differential equations, whose solution
set can in principle be determined. In this fashion we reduce \eqref{SHS} to a Hamiltonian system with
finitely many degrees of freedom which can be treated by well-developed dynamical-systems methods, in
particular normal-form theory. We proceed by perturbing the physical parameters $\beta$, $\alpha$ around fixed reference values
$\beta_0$, $\alpha_0$ and thus introducing bifurcation parameters $\varepsilon_1$,  $\varepsilon_2$.
The Kirchg\"{a}ssner reduction delivers an $\varepsilon$-dependent reduced system
which captures the small-amplitude dynamics for small values of these
parameters; its dimension is the number of purely
imaginary eigenvalues of the corresponding linearised system at $(\varepsilon_1,\varepsilon_2)=(0,0)$.
The reduction procedure is therefore especially helpful in detecting
bifurcations which are associated with a change in the number of
purely imaginary eigenvalues.

\begin{figure}
\hspace{4.5cm}\includegraphics[width=7cm]{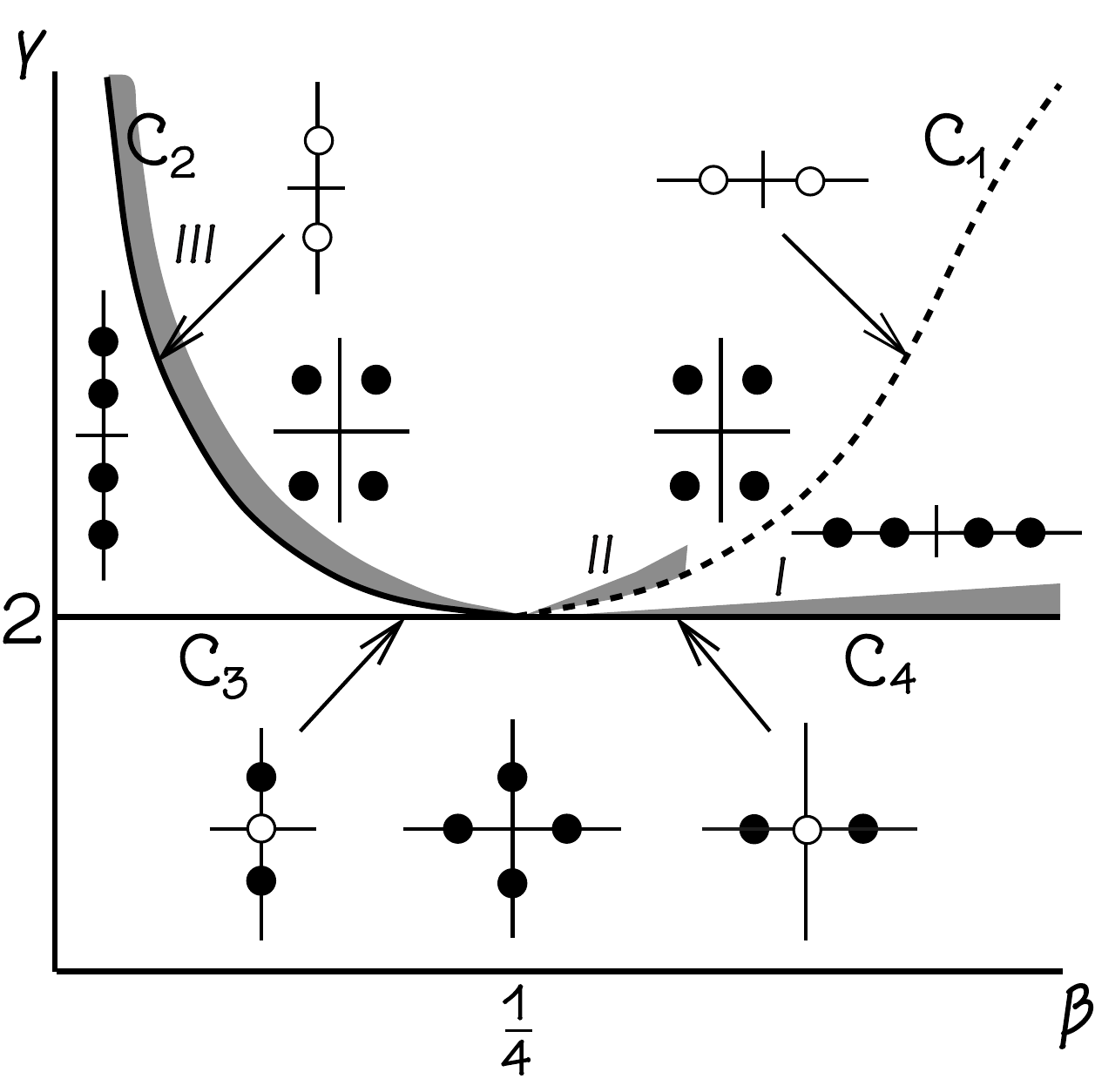}

\vspace*{1mm}
{\it (a) Bifurcation curves in the $(\beta,\gamma)$-plane; the shaded regions indicate the
parameter regimes in which homoclinic bifurcation is detected.

\vspace*{4mm}
\hspace{3.25cm}\includegraphics[width=3.5cm]{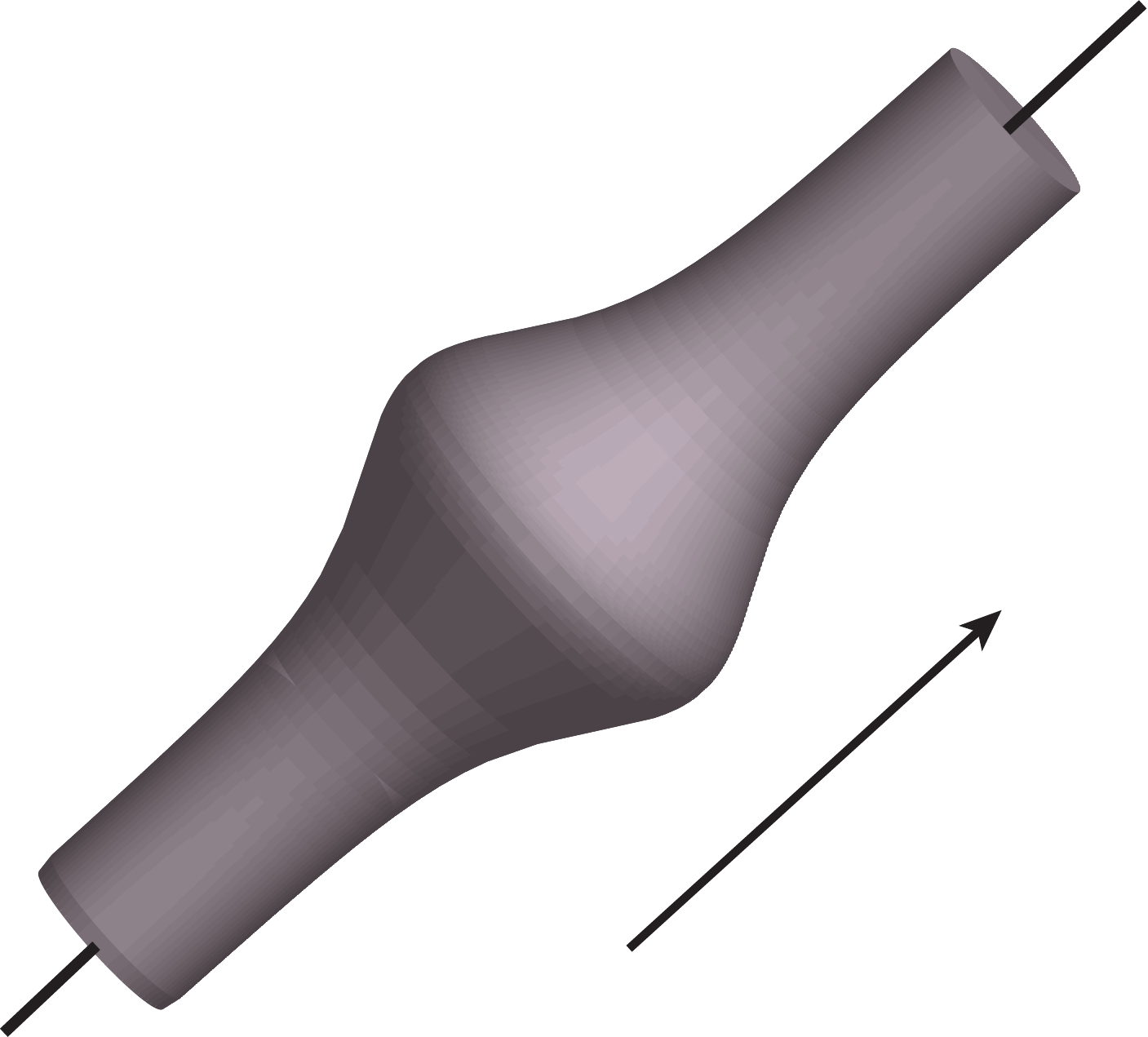}\hspace{2.5cm}\includegraphics[width=3.5cm]{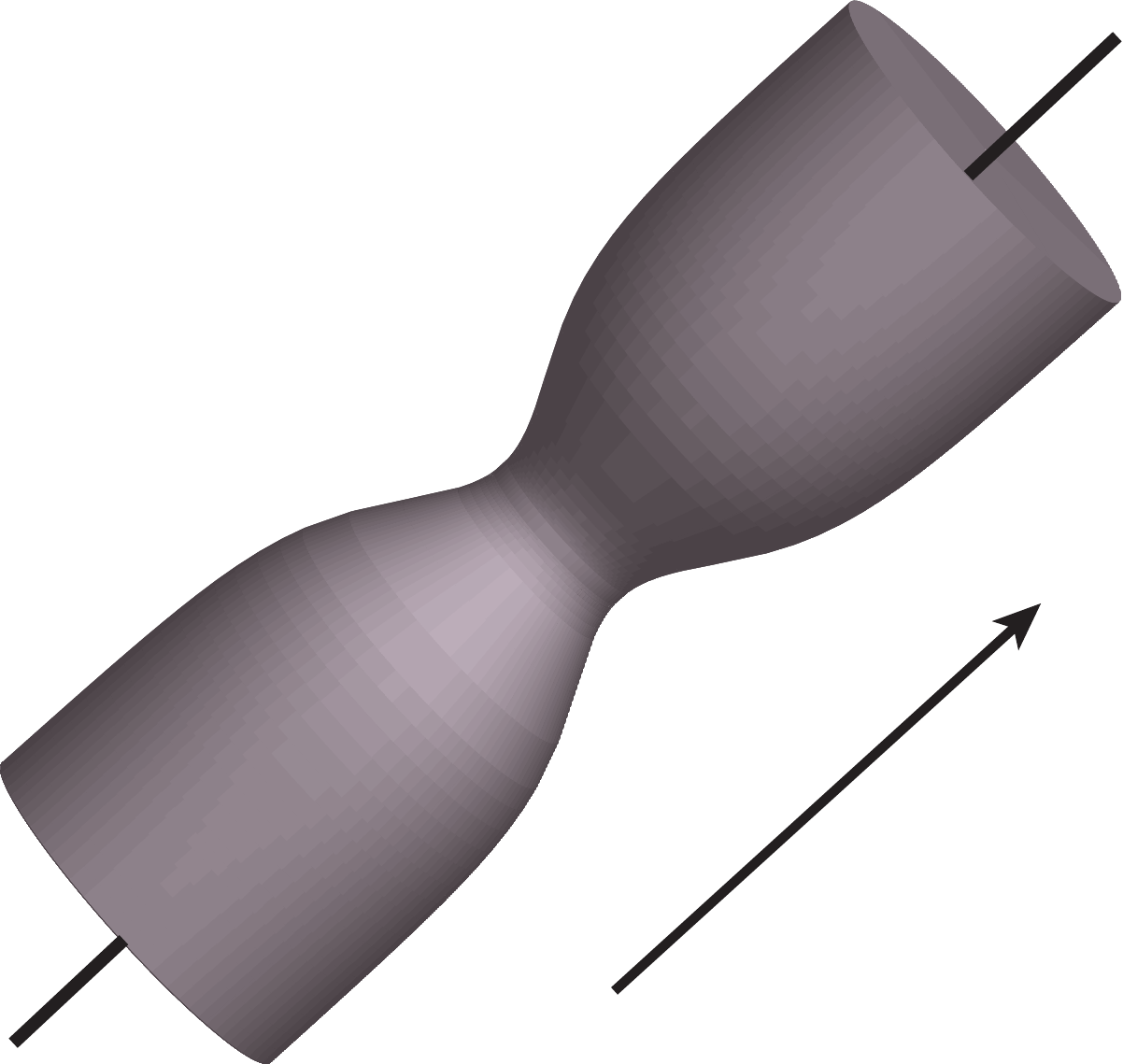}

\vspace*{4mm}
(b) Solitary waves of elevation (left) and depression (right) in region I.

\vspace*{4mm}
\hspace{3.25cm}\includegraphics[width=3.5cm]{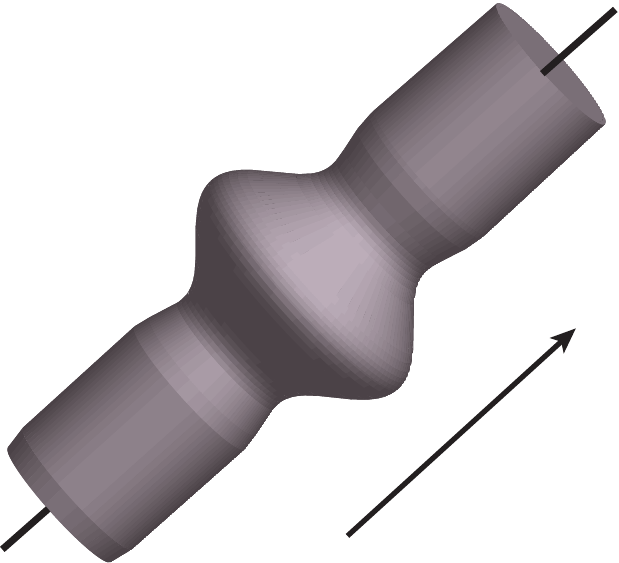}\hspace{2.5cm}\includegraphics[width=3.5cm]{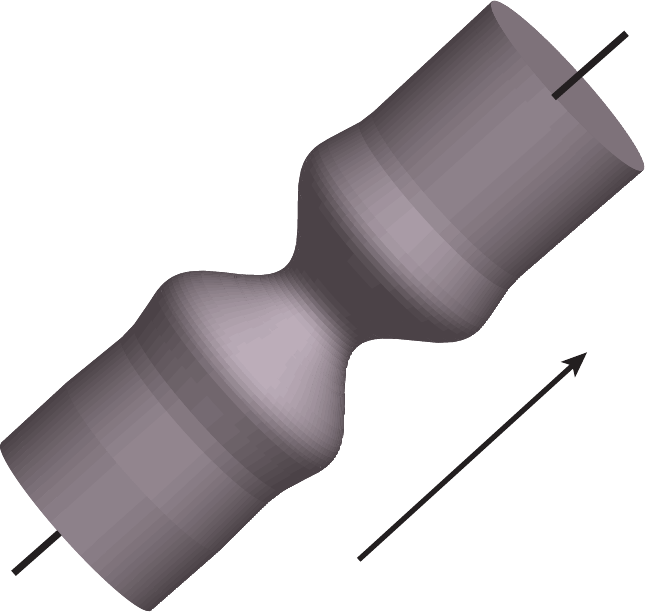}

\vspace*{4mm}
(c) Primary solitary waves of elevation (left) and depression (right) in region II.

\vspace*{4mm}
\hspace{3.25cm}\includegraphics[width=3.5cm]{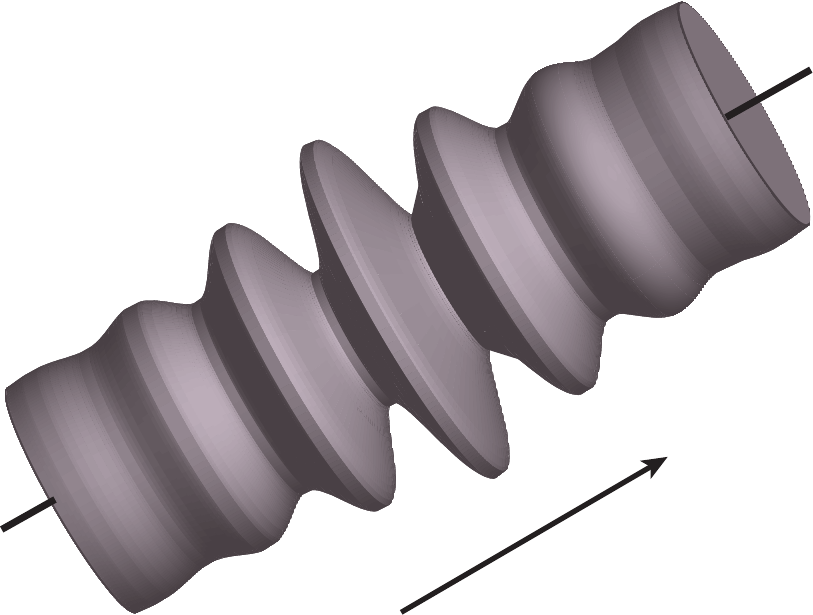}\hspace{2.5cm}\includegraphics[width=3.5cm]{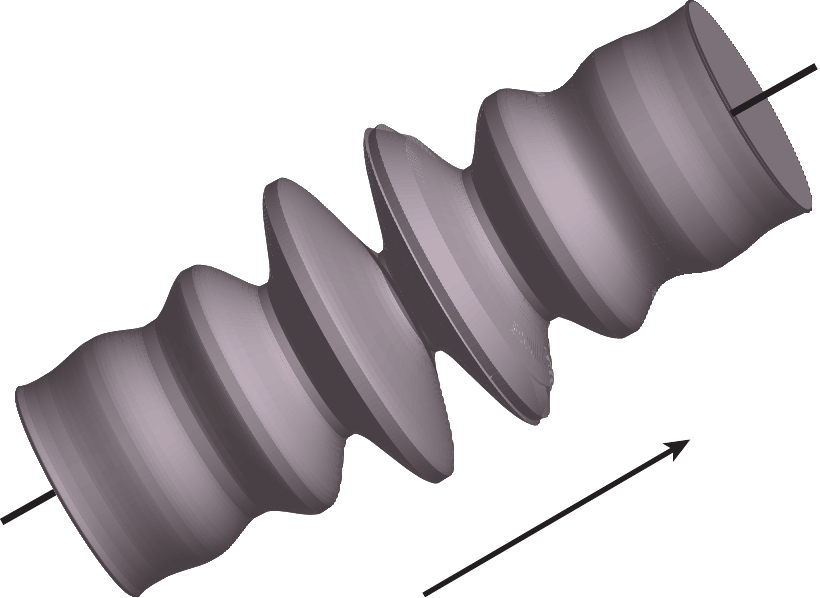}

\vspace*{4mm}
(d) Primary solitary waves of elevation (left) and depression (right) in region III.

\caption{Summary of the basic types of solitary wave whose existence is established in the
present paper by the Kirchg\"{a}ssner reduction.}
\label{summary}}
\end{figure}

Working in the $(\beta,\gamma)$ parameter plane, where $\gamma=\alpha-\beta$,
one finds that there are three critical curves $C_2$, $C_3$, $C_4$
at which the number of purely imaginary
eigenvalues changes (see Figure \ref{summary}(a)),
together with a fourth curve $C_1$ at which the number of real eigenvalues changes.
(In fact $C_3=\{(\beta,2): \beta<\frac{1}{4}\}$,
$C_4=\{(\beta,2): \beta>\frac{1}{4}\}$ and explicit formulae for $C_1$ and $C_2$
are given in Section \ref{Reduction}.) A similar diagram arises in the study of gravity-capillary
travelling water waves (see Iooss \cite{Iooss95}, Groves \& Wahl\'{e}n \cite{GrovesWahlen07} and the references
therein), and there the curves corresponding to $C_1$, $C_2$ and $C_4$ are associated with
\emph{homoclinic bifurcation:} homoclinic solutions of the
reduced Hamiltonian system (corresponding to solitary water waves)
bifurcate from the trivial solution. Figure \ref{summary}(a) illustrates the parameter
regions I, II and III adjacent to $C_1$, $C_2$ and $C_4$
in which the existence of homoclinic solutions is to be expected. In Section \ref{Analysis of reduced systems}
we study these regions using the Kirchg\"{a}ssner reduction; the basic types of solitary wave found there
are sketched in Figures \ref{summary}(b)--(d).

In Section \ref{Region I} we examine region I, choosing $(\beta_0,\gamma_0) \in C_4$, so that
$\alpha_0=2+\beta_0$, and writing $\alpha=\alpha_0+\mu$ with $0 < \mu \ll 1$. According to the Kirchg\"{a}ssner reduction
small-amplitude solitary waves are given by
$$\eta(z)=\tfrac{1}{2}\mu(\beta_0-\tfrac{1}{4})^{1/2}Q\left(\mu^{1/2}(\beta_0-\tfrac{1}{4})^{-1/2}z\right) + O(\mu^{3/2}),$$
where $(Q,P)$ is a homoclinic solution of the reversible Hamiltonian system
\begin{align}
\dot{Q} & = P+O(\mu^{1/2}), \label{KdV 1} \\
\dot{P} & = Q - \check{c}_1 Q^2 + O(\mu^{1/2}) \label{KdV 2}
\end{align}
with $\check{c}_1 := \tfrac{1}{2}(\alpha_0m_1^\prime(1)-6)$.
This system admits a homoclinic solution which corresponds to a monotonically decaying, symmetric solitary
wave of elevation for $\check{c}_1>0$ and depression for $\check{c}_1<0$. For $m_1^\prime(1)$ close to the critical value
$6\alpha_0^{-1}$ we write $m_1^\prime(1)=\alpha_0^{-1}(6+2\mu^{1/2}\check{\kappa})$ with $0<\check{\kappa} \ll 1$
and find that small-amplitude solitary waves are given by
$$\eta(z)=\tfrac{1}{2}\mu^{1/2}(\beta_0-\tfrac{1}{4})^{1/2}Q\left(\mu^{1/2}(\beta_0-\tfrac{1}{4})^{-1/2}z\right) + O(\mu^{3/2}),$$
where $(Q,P)$ is a homoclinic solution of the reversible Hamiltonian system
\begin{align}
\dot{Q} & = P+O(\mu^{1/2}), \label{cubic KdV 1} \\
\dot{P} & = Q - \check{\kappa} Q^2 -\check{d}_1 Q^3+O(\mu^{1/2}) \label{cubic KdV 2}
\end{align}
with $\check{d}_1=\tfrac{1}{6}(12-\alpha_0 m_1^{\prime\prime}(1))$. For $\check{d}_1>0$ this system admits a pair of homoclinic solutions which correspond to monotonically decaying, symmetric solitary
waves; one is a wave of depression, the other a wave of elevation. Note that in the limit $\mu=0$ or
$(\mu,\check{\kappa})=(0,0)$ the variable $Q$ solves a travelling-wave
version of the (generalised) Korteweg-de Vries equation.

In Section \ref{Region II} we apply the Kirchg\"{a}ssner reduction in region II, finding that
small-amplitude solitary waves are given by
$$\eta(z)=\tfrac{1}{2}\mu^4P_1(\mu z) + O(\mu^5),$$
where $(Q,P)$ is a homoclinic solution of the reversible Hamiltonian system
\begin{align*}
\dot{Q}_1 & = -P_1 + \tfrac{2}{3}(1+\delta)P_2 + \tfrac{4}{9}(1+\delta)^2P_1 + 3 c_1P_1^2+ O(\mu), \\
\dot{Q}_2 & = P_2  + \tfrac{2}{3}(1+\delta)P_1 + O(\mu), \\
\dot{P}_1 & = Q_2+ O(\mu), \\
\dot{P}_2 & = Q_1 + \tfrac{2}{3}(1+\delta)Q_2+ O(\mu)
\end{align*}
with $c_1 = 48\sqrt{6}(3 m_1^\prime(1)-8)$; the parameters $0 < \mu, \delta \ll 1$ measure the distance from respectively the point $(\beta_0,\gamma_0)=(\frac{1}{4},2)$
and the curve $C_1$. This system admits a homoclinic solution which corresponds to a solitary
wave of elevation for $c_1>0$ and depression for $c_1<0$; the wave is symmetric with an oscillatory
decaying tail. 
For $m_1^\prime(1)$ close to the critical value
$\frac{8}{3}$ we write $m_1^\prime(1)=\tfrac{1}{3}(8+\frac{1}{144\sqrt{6}}\check{\kappa}\mu^2)$ with $0<\check{\kappa} \ll 1$
and find that small-amplitude solitary waves are given by
$$\eta(z)=\tfrac{1}{2}\mu^2P_1(\mu z) + O(\mu^3),$$
where $(Q,P)$ is a homoclinic solution of the reversible Hamiltonian system
\begin{align*}
\dot{Q}_1 & = -P_1 + \tfrac{2}{3}(1+\delta)P_2 + \tfrac{4}{9}(1+\delta)^2P_1 + \check{\kappa} P_1^2+4d_1P_1^3+ O(\mu), \\
\dot{Q}_2 & = P_2  + \tfrac{2}{3}(1+\delta)P_1 + O(\mu), \\
\dot{P}_1 & = Q_2+ O(\mu), \\
\dot{P}_2 & = Q_1 + \tfrac{2}{3}(1+\delta)Q_2+ O(\mu)
\end{align*}
with $d_1 = 864\left(\frac{1264}{75}-m_1^{\prime\prime}(1)\right)$.
For $d_1>0$ this system admits a
a pair of homoclinic solutions which correspond to symmetric solitary
waves with oscillatory decaying tails; one is a wave of depression, the other a wave of elevation.
Note that in the limit $\mu=0$ or
$(\mu,\check{\kappa})=(0,0)$ the variable $P_1$ solves a travelling-wave
version of the (generalised) Kawahara equation.

It is instructive to interpret the above results for two well-studied magnetic intensities.\\
\\
\emph{(i) The linear magnetisation law}
$$m_1(s)=s.$$
In region I we find that $\check{c}_1<0$ for $\alpha_0<6$
(solitary waves of depression) and $\check{c}_1>0$ for $\alpha_0>6$ (solitary waves of elevation); furthermore
$\check{d}_1=2$, so that both types of waves exist for $\alpha_0$ near $6$. This region has also been studied by
Rannacher \& Engel \cite{RannacherEngel06} using a weakly nonlinear approximation.
In terms of the magnetic Bond number $B=\alpha_0/\beta$ (corresponding to
$B<9$) they derived a Korteweg-de Vries equation equivalent to \eqref{KdV 1}, \eqref{KdV 2} and
found solitary waves of depression for $\frac{3}{2}<B<9$ (that is, $\alpha_0<6$) and of elevation
for $1<B<\frac{3}{2}$ (that is, $\alpha_0>6$), in agreement with our results.
(Continuing their weakly nonlinear analysis to the next order in this
region would lead to a cubic Korteweg-de Vries equation equivalent to \eqref{cubic KdV 1}, \eqref{cubic KdV 2} and the
prediction of both types of waves for $B$ near $\frac{3}{2}$). In region II we find that $c_1 = -240\sqrt{6}$
(solitary waves of depression).\\
\\
\emph{(ii) The Langevin magnetisation law}
$$
m_1(s) = \frac{\coth (\lambda s)-(\lambda s)^{-1}}{\coth \lambda - \lambda^{-1}}, \label{Langevin law}
$$
where $\lambda>0$ is a dimensionless parameter. In Region I we find that $\check{c}_1<0$ for $\alpha_0<6$
and $\alpha_0>6$, $\lambda \in (\lambda^\star(\alpha_0),\infty)$ (solitary waves of depression), while $\check{c}_1<0$ for
$\alpha_0>6$, $\lambda \in (0,\lambda^\star(\alpha_0))$ (solitary waves of elevation), where $\lambda^\star(\alpha_0)$ is the unique solution of
the equation
$$\frac{\lambda^{-1}-\lambda\cosech^2 \lambda}{\coth \lambda -\lambda^{-1}}=6\alpha_0^{-1}$$
(so that $\lambda^\star(6)=0$).
Furthermore $\check{d}_1>0$, so that both types of waves exist for $(\lambda,\alpha_0)$ near $(\lambda^\star,\alpha_0(\lambda^\star))$
(with $\alpha_0(0)=6$).
In region II we find that $c_1<0$ (solitary waves of depression).\pagebreak

In Section \ref{Region III} we turn to region III.
Introducing a bifurcation parameter $\mu$ so that positive
values of $\mu$ correspond to points on the `complex' side of $C_2$,
one obtains the reduced (reversible) Hamiltonian system
$$\dot{A} = \frac{\partial \tilde{H}^\mu}{\partial \bar{B}}, \qquad
\dot{B} = - \frac{\partial \tilde{H}^\mu}{\partial \bar{A}},$$
$$\tilde{H}^\mu=\i s(A\bar{B}-\bar{A}B)+|B|^2+\tilde{H}_\mathrm{NF}^0(|A|^2,\i(A\bar{B}-\bar{A}B),\mu)
+ O(|(A,B)|^2|(\mu,A,B)|^{n_0}),$$
where $\tilde{H}_\mathrm{NF}^0$ is a real polynomial which satisfies
$\tilde{H}_\mathrm{NF}^0=0$; it contains the terms of order
$3$, \ldots, $n_0+1$ in the Taylor expansion of $\tilde{H}^\mu$.
The substitution
$A(z)=\mathrm{e}^{\i s z}a(z)$, $B(z)=\mathrm{e}^{\i s z}b(z)$
converts the `truncated normal form' obtained by neglecting the remainder term into the system
\begin{align*}
\dot{a} & = b + \partial_b\tilde{H}_\mathrm{NF}^0(|a|^2,\i(a\bar{b}-\bar{a}b),\mu) , \\
\dot{b} & = -\partial_{\bar{a}}\tilde{H}_\mathrm{NF}^0(|a|^2,\i(a\bar{b}-\bar{a}b),\mu)
\end{align*}
(which, as evidenced by the scaling $z \mapsto \mu^{1/2}z$, $(a,b) \mapsto (\mu^{1/2}a,\mu b)$, is
at leading order equivalent to the nonlinear Schr\"{o}dinger equation).
Supposing that the coefficients
of certain terms in $\tilde{H}_\mathrm{NF}^0$ have the correct sign, one finds that
the latter system admits
a circle of homoclinic solutions, two of which are real. The corresponding pair of homoclinic solutions
to the original `truncated normal form' are reversible and persist when the remainder terms are
reinstated (see Iooss \& P\'{e}rou\`{e}me \cite{IoossPeroueme93}). They generate
symmetric solitary waves which take the form of periodic wave trains modulated by
exponentially decaying envelopes; one is a wave of depression, the other a wave of elevation.

Each of the basic types of solitary waves in regions II and III is the \emph{primary}
member of an infinite family of \emph{multipulse} solitary waves which resemble
multiple copies of the primary. These waves are generated by corresponding multipulse
homoclinic solutions which make several large excursions away from the origin
in their four-dimensional phase space. A more precise description of the multipulse
waves, together with a discussion of the relevant existence theories (which are
based on variational and dynamical-systems arguments) is given in
Sections \ref{Region II} and \ref{Region III}.

Although the techniques used in the present paper are generalisations of those developed
for the water-wave problem (see Iooss \cite{Iooss95}, Groves \& Wahl\'{e}n \cite{GrovesWahlen07}
and the references therein), we employ different methods to compute the reduced
Hamiltonian systems. The spatial Hamiltonian system \eqref{SHS} is invariant under the transformation
$\hat{\phi} \mapsto \hat{\phi} + c$, $c \in {\mathbb R}$ (`variation of potential base-level'), and the quantity $\int_0^1 r \hat{\zeta} \dr$ is conserved.
In many hydrodynamic problems it is possible to eliminate a symmetry of this kind before applying the Kirchg\"{a}ssner
reduction (see e.g.\ Groves, Lloyd \& Stylianou \cite[\S3.1]{GrovesLloydStylianou16} for an example
in stationary ferrofluids), but here we retain it. It is inherited by the reduced systems: one of the canonical coordinates
is cyclic and its conjugate
is conserved. According to the
classical theory, the next step is to set the conserved variable to zero, solve the resulting decoupled system for the other
variables
and recover the cyclic variable by quadrature; the lower-order system is typically studied using a canonical change
of variables which simplifies its Hamiltonian (a `normal-form' transformation). In the present context it is convenient to use a normal-form
transformation \emph{before} lowering the order of the system since it can be `absorbed' into the changes of variable
associated with the Kirchg\"{a}ssner reduction; this procedure greatly simplifies our later calculations. We
present a general result for this purpose (Theorem \ref{BM theorem}), whose proof is
based upon the method given by Bridges \& Mielke \cite[Theorem 4.3]{BridgesMielke95} and which may also
be helpful in other applications.

\section{The ferrohydrodynamic problem} \label{ferroeqns}
 
We consider an incompressible, inviscid ferrofluid of unit density in the region
$$S_1:=\{0<r<R+\eta(\theta,z,t)\}$$
bounded by the free interface $\{r=R+\eta(\theta,z,t)\}$
and a current-carrying wire at $\{r=0\}$, where $(r,\theta,z)$ are cylindrical polar coordinates.
The fluid is subject to a static magnetic field and
the surrounding region $$S_2=\{r>R+\eta(\theta,z,t)\}$$
is a vacuum (see Figure \ref{Definition sketch}).

We denote the magnetic and induction fields in the fluid and vacuum by respectively
$\mathbf{H}_1, \mathbf{B}_1$ and $\mathbf{H}_2$, $\mathbf{B}_2$, and suppose that the
relationships between them are given by the identities
$$
\mathbf{B}	_1=\mu_{0}(\mathbf{H}_1+\mathbf{M}_1({\bf H}_1)), \qquad
\mathbf{B}_2=\mu_{0}\mathbf{H}_2,
$$
where $\mu_{0}$ is the magnetic permeability of free space and $\bf{M}_1$ is the (prescribed)
magnetic intensity of the ferrofluid. We suppose that
$$\mathbf{M}_1({\bf H}_1)=m_1(|{\bf H}_1|)\frac{{\bf H}_1}{|{\bf H}_1|}$$
where $m_1$ is a (prescribed) nonnegative function, so that in particular ${\bf M}_1$ and ${\bf H}_1$ are collinear.
According to Maxwell's equations the magnetic and induction fields are respectively irrotational
and solenoidal, and introducing  magnetic potential functions $\psi_1$, $\psi_2$ with ${\bf H}_1=-\nabla \psi_1$,\linebreak
${\bf H}_2=-\nabla \psi_2$, we therefore find that
\begin{align*}
\nabla\cdot(\mu(|\nabla \psi_1|)\nabla \psi_1) & =0 \qquad \mbox{in $S_1$}, \\
\Delta \psi_2 & =0 \qquad\mbox{in $S_2$},
\end{align*}
in which
$$\mu(s)=1+\frac{m_1(s)}{s}$$
is the magnetic permeability of the ferrofluid relative to that of free space.
We suppose that the ferrofluid flow is irrotational, so that its velocity field ${\bf v}$ is
the gradient of a scalar velocity potential $\phi$. The Euler equation for the ferrofluid is given by
$${\bf v}_t + ({\bf v}.\nabla){\bf v} = - \nabla p^\star + \mu_{0}
			\left(\mathbf{M}_1\cdot
			\nabla\right) {\mathbf{H}_1}
$$
(Rosensweig \cite[\S5.1]{Rosensweig}), where $p^\star$ is its composite pressure, and
the calculations
$$\left(\bf{M}_1\cdot\nabla\right){\mathbf{H}_1}=|{\bf M}_1|\nabla (|{\bf H}_1|) = \nabla \left( \int_0^{|{\bf H}_1|} m_1(t)\dt\right), \qquad
({\bf v}.\nabla){\bf v} = \nabla \left( \frac{1}{2}|{\bf v}|^2\right)$$
show that this equation is equivalent to
\begin{equation}
\phi_t + \frac{1}{2} |\nabla \phi|^2 - \mu_0\int_0^{|{\bf H}_1|} m_1(t)\dt+p^{\star} =c_{0}, \label{Euler}
\end{equation}
where $c_0$ is a constant.

Next we turn to the boundary conditions at $\{r=R+\eta(\theta,z,t)\}$. The \emph{magnetic boundary conditions} are
$$
{\bf H}_1\cdot \mathbf{t}=\bf{H}_2\cdot \mathbf{t},
\qquad
\bf{B}_1\cdot \mathbf{n}=\bf{B}_2\cdot \mathbf{n},		
$$
where ${\bf t}$ and ${\bf n}$ denote tangent and normal vectors to the free surface; it follows that
$$
\psi_2-\psi_1\Big|_{r=R+\eta(\theta,z,t)}=0,\quad \psi_{2n}-\mu(|\nabla\psi_1|) \psi_{1n}\Big|_{r=R+\eta(\theta,z,t)} =0.
$$
The \emph{(hydro-)dynamical boundary condition} is given by 
$$p^{\star}+\frac{\mu_{0}}{2}({\mathbf{M}_1\cdot \mathbf n})^{2}= 2\sigma\kappa,$$
(Rosensweig \cite[\S5.2]{Rosensweig}),
in which $\sigma>0$ is the coefficient of surface tension and
\small
$$2\kappa = \frac{-2\eta_\theta^2 - (R+\eta)^2(1+\eta_z^2)+(R+\eta)^3\eta_{zz}+(R+\eta)\eta_\theta^2\eta_{zz}
-2(R+\eta)\eta_\theta\eta_z\eta_{\theta z}+(R+\eta)(1+\eta_z^2)\eta_{\theta\theta}}
{((R+\eta)(1+\eta_z^2)+\eta_\theta^2)^{3/2}}$$
\normalsize
is the mean curvature of the interface; using \eqref{Euler}, we find that
$$
\phi_t + \frac{1}{2} |\nabla \phi|^2 - \mu_0\nu(|\nabla \psi_1|) + 2\sigma \kappa
-\frac{\mu_{0}}{2}(\mu(|\nabla \psi_1|)-1)^2 \Bigg|_{r=R+\eta(\theta,z,t)} = c_0,
$$
where
$$\nu(s) = \int_0^s m_1(t)\dt.$$
Finally,  the \emph{(hydro-)kinematic boundary condition} is
$$(\partial_t + {\bf v}.\nabla)(r-R-\eta(\theta,z,t))=0,$$
that is
$$-\eta_t + \phi_r - \frac{1}{r^2}\phi_\theta\eta_\theta - \phi_z \eta_z=0\Bigg|_{r=R+\eta(\theta,z,t)}.$$
The relevant conditions at $r=0$ and in the far field are ${\bf v}.{\bf e}_r$, ${\bf B}_1.{\bf e}_r
\rightarrow 0$ as $r \rightarrow 0$, so that $\phi_r$, $\psi_{1r} \rightarrow 0$
as $r \rightarrow 0$, and ${\bf B}_2.{\bf e}_r \rightarrow 0$ as $r \rightarrow \infty$,
so that $\psi_{2r} \rightarrow 0$ as $r \rightarrow \infty$.

The constant $c_0$ is selected so that
$${\bf H}_1=\dfrac{J}{2\pi r}{\bf e}_\theta, \qquad{\bf H}_2=\dfrac{J}{2\pi r}{\bf e}_\theta, \qquad {\bf v}={\bf 0}, \qquad \eta=0$$
(that is $\psi_1=\psi_2=-J\theta/2\pi$, $\phi=0$, $\eta=0$)
is a solution to the above equations
(corresponding to a
uniform magnetic field and a circular cylindrical jet with radius $R$); we therefore set
$c_0 = -\mu_0\nu(J/2\pi r)+\sigma/R$.
Seeking axisymmetric waves for which $\eta$ and $\phi$ are independent of $\theta$, one finds that
$\psi_1=\psi_2=-J\theta/2\pi$, so that the hydrodynamic problem decouples from the magnetic problem
and is given by
\begin{align*}
&\phi_{rr}+\frac{1}{r}\phi_r + \phi_{zz} =0, \qquad 0<r<R+\eta(z,t), \\
& \phi_ r = 0, \hspace{1.29in} r = 0
\end{align*}
and
\begin{align*}
& -\eta_t+\phi_r-\phi_z\eta_z = 0, \\
& \phi_t+\frac{1}{2}(\phi_r^2+\phi_z^2) \\
& \hspace{1cm}\mbox{}-\mu_0\nu\left(\frac{J}{2\pi(R+\eta)}\right)
+\mu_0\nu\left(\frac{J}{2\pi R}\right)+\frac{\sigma}{(R+\eta)(1+\eta_z^2)^{1/2}} - \frac{\sigma \eta_{zz}}{(1+\eta_z^2)^{3/2}}-\frac{\sigma}{R}=0
\end{align*}
for $r=R+\eta(z,t)$.

The next step is to seek travelling wave solutions for which $\eta$ and $\phi$ depend upon $z$ and $t$
only through the combination $z-ct$, and to introduce dimensionless variables
\[
 (\hat z,\hat r):=\frac{1}{R}(z-ct,r),\quad \hat \phi:= \frac{1}{cR}\phi,\quad
\hat\eta:= \frac{1}{R}\eta.
\]
and functions
\[
\hat{m}_1(s):=\frac{2\pi R}{J \chi}m_1\left(\frac{J}{2\pi R}s\right), \qquad
\hat\nu(s):= \frac{4\pi^2 R^2}{J^2\chi}\nu\left(\frac{J}{2\pi R}s\right),
\]
where
\[
\chi=\frac{2\pi R}{J}m_1\left(\frac{J}{2\pi R}\right)
\]
(note that $\hat{m}(1)=\hat{\nu}^\prime(1)=1$). Dropping the hats for notational simplicity, we find that
\begin{align}
&\phi_{rr}+\frac{1}{r}\phi_r + \phi_{zz} =0, \qquad 0<r<1+\eta(z,t), \label{Unflattened 1} \\
& \phi_ r = 0, \hspace{1.29in} r = 0 \label{Unflattened 2}
\end{align}
and
\begin{align}
& \eta_z+\phi_r-\phi_z\eta_z = 0, \label{Unflattened 3}\\
& -\phi_z+\frac{1}{2}(\phi_r^2+\phi_z^2)-\alpha \frac{T^\prime(\eta)}{1+\eta}
+\beta\left(\frac{1}{(1+\eta)(1+\eta_z^2)^{1/2}} - \frac{\eta_{zz}}{(1+\eta_z^2)^{3/2}}-1\right)=0 \label{Unflattened 4}
\end{align}
for $r=1+\eta(z,t)$, where 
$$T(\eta)=\int_0^\eta \left(\nu\left(\frac{1}{1+s}\right)-\nu(1)\right)(1+s)\ds$$
and
\begin{equation}
\alpha=\frac{\mu_0 J^2\chi}{4\pi^2R^2c^2}, \qquad \beta = \frac{\sigma}{c^2R}
\label{Defns of alpha and beta}
\end{equation}
are dimensionless parameters. \emph{Solitary waves} are nontrivial solutions of \eqref{Unflattened 1}--\eqref{Unflattened 4}
with $\eta(z)$, $\phi(r,z) \rightarrow 0$ as $z \rightarrow \pm \infty$.
Finally, note that equations \eqref{Unflattened 1}, \eqref{Unflattened 3} and \eqref{Unflattened 4} follow from the formal variational principle
$$\delta \int \left\{
\int_0^{1+\eta}\!\!\left(\frac{1}{2}r\phi_r^2+\frac{1}{2}r\phi_z^2-r\phi_z\right)\!\!\dr-\alpha T(\eta)+\beta(1+\eta)(1+\eta_z^2)^{1/2}-\frac{1}{2}\beta(1+\eta)^2\right\}\dz=0,$$
where the variations are taken with respect to $\eta$ and $\phi$.

\section{Spatial dynamics} \label{Spatial dynamics}

\subsection{Formulation as a spatial Hamiltonian system} \label{SD - derivation}
The first step is to use the `flattening' transformation
$$
\hat{r}=\frac{r}{1+\eta}
$$
to map the variable domain $\{0<r<1+\eta\}$ into a fixed strip $(0,1) \times {\mathbb R}$
and the free interface $\{r=1+\eta(z)\}$ into $\{\hat{r}=1\}$. Dropping the hat for
notational simplicity, we find that the corresponding `flattened' variable
$$\hat{\phi}(\hat{r},z)=\phi(r,z)$$
satisfies the equations
\begin{equation}
\frac{(r\phi_r)_r}{(1+\eta)^2}+r\left(\phi_z-\frac{r\eta_z\phi_r}{1+\eta}\right)_z
-\frac{r^2\eta_z}{1+\eta}\left(\phi_z-\frac{r\eta_z\phi_r}{1+\eta}\right)_r = 0, \qquad 0<r<1
\label{Flattened 1}
\end{equation}
with boundary conditions
\begin{equation}
\phi_r|_{r=0}=0 \label{Flattened 2}
\end{equation}
and
\begin{align}
& \eta_z+\frac{\phi_r}{1+\eta}-\left(\phi_z-\frac{r\eta_z\phi_r}{1+\eta}\right)\eta_z\Bigg|_{r=1} = 0, \label{Flattened 3}\\
& -\left(\phi_z-\frac{r\eta_z\phi_r}{1+\eta}\right)+\frac{1}{2(1+\eta)^2}\phi_r^2 + \frac{1}{2}\left(\phi_z-\frac{r\eta_z\phi_r}{1+\eta}\right)^2\Bigg|_{r=1} \nonumber \\
& \hspace{2cm}\mbox{}-\alpha \frac{T^\prime(\eta)}{1+\eta}
+\beta\left(\frac{1}{(1+\eta)(1+\eta_z^2)^{1/2}} - \frac{\eta_{zz}}{(1+\eta_z^2)^{3/2}}-1\right)=0. \label{Flattened 4}
\end{align}
Observe that equations \eqref{Flattened 1}, \eqref{Flattened 3} and \eqref{Flattened 4} follow from the new variational principle
$\delta \LL=0$, where
\begin{align}
\LL(\eta,\phi)&:=\int
\Bigg\{\int_0^1 \left\{\frac{1}{2}\left(r\phi_r^2 + \left(\phi_z-\frac{r\eta_z\phi_r}{1+\eta}\right)^2(1+\eta)^2r \right)
-\left(\phi_z-\frac{r\eta_z\phi_r}{1+\eta}\right)(1+\eta)^2r\right\}\dr \nonumber \\
&\hspace{0.75in}\mbox{}-\alpha T(\eta)+\beta(1+\eta)(1+\eta_z^2)^{1/2}-\frac{1}{2}\beta(1+\eta)^2\Bigg\}\dz
\label{Definition of cal L}
\end{align}
 and the variations are taken in $\eta$ and $\phi$ (the functional $\LL$ is obtained from the variational functional
 for \eqref{Unflattened 1}, \eqref{Unflattened 3} and \eqref{Unflattened 4} by `flattening').

We exploit this variational principle by regarding $\LL$ as an action functional
of the form
$$\LL = \int L(\eta, \phi, \eta_z, \phi_z) \dz,$$
in which $L$ is the integrand on the right-hand side of equation \eqref{Definition of cal L},
and deriving a canonical Hamiltonian formulation of \eqref{Flattened 1}--\eqref{Flattened 4} by means of the Legendre transform.
To this end, let us introduce new variables $\omega$ and $\xi$ by the formulae
\begin{align*}
\omega &= \frac{\delta L}{\delta \eta_z}
=\int_0^1 \left\{-\left(\phi_z - \frac{r\eta_z\phi_r}{1+\eta}\right)(1+\eta)r^2\phi_r + (1+\eta)r^2\phi_r \right\}\dr
+ \beta \frac{(1+\eta)\eta_z}{(1+\eta_z^2)^{1/2}}, \\
\xi &= \frac{\delta L}{\delta \phi_z} = \left(\phi_z-\frac{r\eta_z\phi_r}{1+\eta}\right)(1+\eta)^2-(1+\eta)^2
\end{align*}
and define the Hamiltonian function by
\begin{align}
H&(\eta,\omega,\phi,\xi) \nonumber\\
& = \eta_z \omega + \int_0^1 r \phi_z \xi \dr - L(\eta, \phi, \eta_z, \phi_z) \nonumber \\
& =\!\! \int_0^1\!\!\left\{ \frac{1}{2}\!\left(\!\frac{\xi}{(1+\eta)^2}+1\!\right)^{\!\!2}\!\!(1+\eta)^2r - \frac{1}{2}r\phi_r^2\right\}\!\!\dr
\!+\!\alpha T(\eta)\! -\!(1+\eta)\sqrt{\beta^2-W^2}+\frac{1}{2}\beta (1+\eta)^2, \label{Hamiltonian}
\end{align}
in which
$$
W=\frac{1}{1+\eta}\left(\omega+\frac{1}{1+\eta}\int_0^1 r^2 \phi_r \xi\dr\right).
$$

Writing $(\beta,\alpha)=(\beta_0+\varepsilon_1,\alpha_0+\varepsilon_2)$, where $(\beta_0,\alpha_0)$ are fixed,
and $\xi=\zeta-1$ (since $(\eta,\omega,\phi,\xi)=(0,0,-1,0)$ is the `trivial' solution of Hamilton's equations),
we find that Hamilton's equations are given explicitly by
\begin{align}
\eta_z = \frac{\delta H^\varepsilon}{\delta \omega} &= \frac{W}{\sqrt{(\beta_0+\varepsilon_1)^2-W^2}}, \label{HE 1}\\
\omega_z = -\frac{\delta H^\varepsilon}{\delta \eta} &
=\int_0^1 \left\{\left(\frac{(\zeta-1)^2}{(1+\eta)^4}-1\right)(1+\eta)r+ \frac{Wr^2 \phi_r (\zeta-1)}{(1+\eta)^2\sqrt{(\beta_0+\varepsilon_1)^2-W^2}}\right\}\dr
\nonumber \\
& \qquad \mbox{}-(\alpha_0+\varepsilon_2) T^\prime(\eta) + \frac{(\beta_0+\varepsilon_1)^2}{\sqrt{(\beta_0+\varepsilon_1)^2-W^2}}-(\beta_0+\varepsilon_1)(1+\eta), \label{HE 2}\\
\phi_z = \frac{\delta H^\varepsilon}{\delta \zeta} &= \left(\frac{(\zeta-1)}{(1+\eta)^2}+1\right)+\frac{W}{\sqrt{(\beta_0+\varepsilon_1)^2-W^2}}\frac{r\phi_r}{1+\eta}, \label{HE 3}\\
\zeta_z = -\frac{\delta H^\varepsilon}{\delta \phi} &=-\frac{1}{r}(r\phi_r)_r+\frac{W}{\sqrt{(\beta_0+\varepsilon_1)^2-W^2}}\frac{1}{r}\frac{(r^2(\zeta-1))_r}{1+\eta},\label{HE 4}
\end{align}
where the superscript denotes the dependence upon $\varepsilon$,
with boundary condition
\begin{equation}
r\phi_r-\frac{W}{\sqrt{(\beta_0+\varepsilon_1)^2-W^2}}\frac{r^2(\zeta-1)}{1+\eta}\Bigg|_{r=1}=0, \label{BC 2}
\end{equation}
the second of which arises from the integration by parts necessary to compute \eqref{HE 4}
Note that our equations are reversible, that is invariant under the transformation
$(\eta,\omega,\phi,\zeta)(z) \mapsto S(\eta,\omega,\phi,\zeta)(-z)$, where
the \emph{reverser} is
defined by $S(\eta,\omega,\phi,\zeta)=(\eta,-\omega,-\phi,\zeta)$.

To make this construction rigorous  we recall the differential-geometric definitions of a Hamiltonian system
and Hamilton's equations for its associated vector field.

\begin{definition} \label{HS definitions}

A \underline{Hamiltonian system} consists of
a triple $(M,\Omega,H)$, where $M$ is a manifold, $\Omega: TM \times TM \rightarrow {\mathbb R}$
is a closed, weakly nondegenerate bilinear form (the \underline{symplectic $2$-form}) and the
\underline{Hamiltonian}
$H: N \rightarrow {\mathbb R}$ is a smooth function on a manifold domain $N$ of $M$ (that is, a
manifold $N$ which is smoothly embedded in $M$ and has the property that $TN|_{n}$ is densely embedded
in $TM|_{n}$ for each $n \in N$).

Its \underline{Hamiltonian vector field} $v_H$ with domain $\DD(v_H) \subseteq N$ is defined
as follows. The point $n \in N$ belongs to $\DD(v_H)$ with $v_H|_{n} := w\in TM|_{n}$ if and only if
$$\Omega|_{n}(w, v) = {\bf d}H|_{n}(v)$$
for all tangent vectors $v \in TM|_{n}$ (by construction ${\bf d}H|_{n} \in T^\ast N|_{n}$ admits a unique extension
${\bf d}H|_{n} \in T^\star M|_{n}$). \underline{Hamilton's equations} for $(M,\Omega,H)$ are the differential equations
$$\dot{u} = v_H|_u$$
which determine the trajectories $u \in C^1({\mathbb R},M) \cap C({\mathbb R},N)$ of its Hamiltonian vector field.
\end{definition}

Definition \ref{HS definitions} applies to the above formulation. Note that the
identity mapping is\linebreak
(up to the scaling factor $\sqrt{2\pi}$) an isometry $\check{L}^2(B_1(0)) \rightarrow L_r^2(0,1)$, $\check{H}^1(B_1(0)) \rightarrow H_r^1(0,1)$
and $\check{H}^2(B_1(0)) \rightarrow \{\phi \in H_r^2(0,1): \phi_r \in L^2_{r^{-1}}(0,1)\}$,
where $B_1(0)$ is the unit ball in ${\mathbb R}^2$ and
$\check{H}^s(B_1(0))$ denotes the closed subspace of $H^s(B_1(0))$ consisting of axisymmetric functions
(see Bernardi, Dauge \& Maday \cite[Theorem II.2.1]{BernardiDaugeMaday}).
We therefore let $M$ be a neighbourhood of the origin in
$$X:=\{(\eta,\omega,\phi,\zeta) \in {\mathbb R} \times {\mathbb R} \times H_r^1(0,1) \times L_r^2(0,1)\}$$
and $N=Y \cap M$ with
$$Y:=\{(\eta,\omega,\phi,\zeta) \in {\mathbb R} \times {\mathbb R} \times H_r^2(0,1) \times H_r^1(0,1): \phi_r \in L^2_{r^{-1}}(0,1)\},$$
so that elements $(\eta,\omega,\phi,\zeta) \in Y$ satisfy $\phi_r|_{r=0}=0$
(see Bernardi, Dauge \& Maday \cite[Remark II.1.1]{BernardiDaugeMaday}).
We consider
values of $(\varepsilon_1,\varepsilon_2)$ in a neighbourhood $\Lambda$ of the origin in ${\mathbb R}^2$ and
choose $M$ and $\Lambda$ small enough so that
$$|\varepsilon_1| < \frac{\beta_0}{4}, \quad \eta > -\frac{1}{2}>-1, \quad |W|<\frac{\beta_0}{2}<\beta_0+\varepsilon_1.$$

The formula
$$
\Omega((\eta_1,\omega_1,\phi_1,\zeta_1), (\eta_2,\omega_2,\phi_2,\zeta_2))
=\omega_2\eta_1-\eta_2\omega_1 + \int_0^1 r(\zeta_2\phi_1-\phi_2\zeta_1) \dr
$$
defines a weakly nondegenerate bilinear form $M \times M \rightarrow {\mathbb R}$ and hence a constant
symplectic $2$-form $TM  \times TM \rightarrow {\mathbb R}$ (its closure follows from the fact that it is constant),
and the function $H^\varepsilon$ given by \eqref{Hamiltonian} belongs to
$C^\infty(N,{\mathbb R})$, so that the triple $(M,\Omega,H)$ is a Hamiltonian system.
Applying the criterion in the definition, one finds that
$$\DD(v_{H^\varepsilon})=
\left\{(\eta,\omega,\phi,\zeta) \in N:r \phi_r-\frac{W}{\sqrt{(\beta_0+\varepsilon_1)^2-W^2}}\frac{r^2(\zeta-1)}{1+\eta}\Bigg|_{r=1}=0\right\}$$
and that Hamilton's equations are given explicitly by \eqref{HE 1}--\eqref{HE 4}.

It remains to confirm the relationship between a solution to Hamilton's equations for\linebreak
$(M,\Omega,H^\varepsilon)$
and a solution to the `flattened' hydrodynamic problem \eqref{Flattened 1}--\eqref{Flattened 4}. 
Suppose that $(\eta,\omega,\phi,\zeta)$ is a smooth
solution of Hamilton's equations. An explicit calculation shows that the variables
$\tilde{\eta}$, $\tilde{\phi}$ given by $\tilde{\eta}(z)=\eta(z)$,
$\tilde{\phi}(r,z) = \phi(z)(r)$ solve 
\eqref{Flattened 1}--\eqref{Flattened 4} (see Groves \& Toland \cite[pp.\ 212-214]{GrovesToland97} for a
discussion of this procedure in the context of water waves).

\section{Centre-manifold reduction} \label{Reduction}

Our strategy in finding solutions to Hamilton's equations
\eqref{HE 1}--\eqref{HE 4} for $(M,\Omega,H^\varepsilon)$ consists in
applying a reduction principle which asserts that $(M,\Omega,H^\varepsilon)$ is locally
equivalent to a finite-dimensional Hamiltonian system. The key result is the following theorem,
which is a parametrised, Hamiltonian version of a reduction principle for quasilinear evolutionary
equations presented by Mielke \cite[Theorem 4.1]{Mielke88a} (see Buffoni,
Groves \& Toland \cite[Theorem 4.1]{BuffoniGrovesToland96}).

\begin{theorem}\label{Mielke's theorem}
Consider the differential equation
\begin{equation}
\dot{u} = \LL u + \NN(u;\lambda), \label{feq}
\end{equation}
which represents Hamilton's equations for the reversible Hamiltonian system 
$(M,\Omega^\lambda,H^\lambda)$.
Here $u$ belongs to a Hilbert space $\XX$, $\lambda \in {\mathbb R}^\ell$ is
a parameter and $\LL: \DD(\LL) \subset \XX \rightarrow \XX$ is a densely 
defined, closed
linear operator. Regarding $\DD(\LL)$ as a Hilbert space equipped with 
the graph norm,
suppose that $0$ is an equilibrium point of \eqref{feq} when $\lambda=0$ and that
\begin{enumerate}
\item[(H1)] The part of the spectrum $\sigma(\LL)$ of $\LL$ which lies on 
the imaginary
axis consists of a finite number of eigenvalues of finite multiplicity
and is separated from the
rest of $\sigma(\LL)$ in the sense of Kato, so that $\XX$ admits the
decomposition $\XX=\XX_1\oplus \XX_2$, where $\XX_1=\PP(\XX)$, 
$\XX_2=(I-\PP)(\XX)$ and
$\PP$ is the spectral projection corresponding the purely imaginary part
of $\sigma(\LL)$.
\item[(H2)] The operator $\LL_2=\LL|_{\XX_2}$ satisfies the estimate
$$\|(\LL_2-\i s I)^{-1}\|_{\XX_2 \rightarrow \XX_2} \leq \frac{C}{1+|s|}, \qquad s \in {\mathbb R},$$
for some constant $C$ that is independent of $s$.
\item[(H3)] There exists a natural number $k$ and neighbourhoods 
$\Lambda \subset
{\mathbb R}^\ell$ of $0$ and $U \subset \DD(\LL)$ of $0$ such that $\NN$ is $(k+1)$
times continuously differentiable on $U \times \Lambda$, its 
derivatives are bounded
and uniformly continuous on $U \times \Lambda$ and $\NN(0,0)=0$,
$\mathrm{d}_1\NN[0,0]=0$.
\end{enumerate}
Under these hypotheses there exist neighbourhoods $\tilde{\Lambda} \subset
\Lambda$ of $0$ and $\tilde{U}_1 \subset U \cap \XX_1$, $\tilde{U}_2 
\subset U \cap \XX_2$
of $0$ and a reduction function
$r:\tilde{U}_1\times\tilde\Lambda\to \tilde{U}_2$ with the following properties.
The reduction function $r$ is $k$ times continuously differentiable
on $\tilde{U}_1\times\tilde\Lambda$, its derivatives are bounded and
uniformly continuous on $\tilde{U}_1\times\tilde\Lambda$ and
$r(0;0)=0$, $\mathrm{d}_1r[0;0]=0$. The graph
$\tilde{M}^\lambda=\{u_1+r(u_1;\lambda) \in \XX_1 \oplus \XX_2: u_1 \in \tilde{U}_1\}$
is a Hamiltonian centre manifold for \eqref{feq}, so that
\begin{list}{(\roman{count})}{\usecounter{count}}
\item
$\tilde{M}^\lambda$ is a locally invariant manifold of \eqref{feq}:
through every point in $\tilde{M}^\lambda$  there passes a unique solution
of \eqref{feq} that remains on $\tilde{M}^\lambda$ as long as it
remains in $\tilde{U}_1\times \tilde{U}_2$.
\item
Every small bounded solution $u(x)$, $x\in{\mathbb R}$ of
\eqref{feq} that satisfies $(u_1(x),u_2(x))\in \tilde{U}_1\times \tilde{U}_2$
lies completely in $\tilde{M}^\lambda$.
\item
Every solution $u_1:(x_1,x_2) \to \tilde{U}_1$ of the reduced equation
\begin{equation}
\dot{u}_1=\LL u_1+\PP\NN(u_1+r(u_1;\lambda);\lambda) \label{req}
\end{equation}
generates a solution
\begin{equation}
u(x)=u_1(x)+r(u_1(x);\lambda) \label{waves}
\end{equation}
of the full equation \eqref{feq}.
\item
$\tilde{M}^\lambda$ is a symplectic submanifold of $M$ and the flow determined
by the Hamiltonian system $(\tilde{M}^\lambda,\tilde{\Omega}^\lambda,
\tilde{H}^\lambda)$,
where the tilde denotes  restriction to $\tilde{M}^\lambda$, coincides
with the flow on $\tilde{M}^\lambda$ determined by
$(M,\Omega^{\lambda},H^\lambda)$.
The reduced equation \eqref{req} is reversible and represents Hamilton's equations for
$(\tilde{M}^\lambda,\tilde{\Omega}^\lambda, \tilde{H}^\lambda)$.
\end{list}
\end{theorem}

Mielke's theorem cannot be applied directly to \eqref{HE 1}--\eqref{HE 4}
because of the nonlinear boundary condition \eqref{BC 2} in the domain of the Hamiltonian vector field
$v_{H^\varepsilon}$ (the right-hand sides of \eqref{HE 1}--\eqref{HE 4} define a smooth mapping
$g^\varepsilon: Y \rightarrow X$ with
$v_{H^\varepsilon}|_u=g^\varepsilon(u)$ for any $u \in \DD(v_{H^\varepsilon})$).
We overcome this
difficulty using the change of variable $G^\varepsilon: (\eta,\omega,\phi,\zeta) \mapsto (\eta,\hat{\omega},\hat{\phi},\zeta)$, where
\begin{equation}
\hat{\omega} = \int_0^1 r^2\phi_r \dr, \qquad
\hat{\phi}=\phi-\frac{W}{\sqrt{(\beta_0+\varepsilon_1)^2-W^2}}\frac{1}{1+\eta}\int_0^r s(\zeta-1)\ds,
\label{Formula for G}
\end{equation}
which transforms the boundary condition in $\DD(v_{H^\varepsilon})$ into
$$r\hat{\phi}_r|_{r=1}=0.$$

\begin{lemma}
For each $\varepsilon \in \Lambda$ the mapping $G^{\varepsilon}$ is a smooth
diffeomorphism from the neighbourhood $M$ of the origin in $X$ onto a neighbourhood
$\hat{M}$ of the origin in $X$, and from $N=M \cap Y$ onto $\hat{N}=\hat{M} \cap Y$.
The diffeomorphisms and their inverses depend smoothly upon $\varepsilon \in \Lambda$.
\end{lemma}
{\bf Proof.} 
These results follow from the explicit formulae \eqref{Formula for G} and
\begin{align*}
\omega&=\frac{(\beta_0+\varepsilon_1)\Gamma(1+\eta)}{\sqrt{1+\Gamma^2}}-\frac{1}{1+\eta}\int_0^1 r^2\hat{\phi}_r(\zeta-1)\dr - \frac{\Gamma}{(1+\eta)^2}\int_0^1 r^3(\zeta-1)^2\dr, \\
\phi & = \hat{\phi} + \frac{\Gamma}{1+\eta}\int_0^r s(\zeta-1)\ds,
\end{align*}
where
$$\Gamma=(1+\eta)\left(\hat{\omega} -\int_0^1 r^2 \hat{\phi}_r \dr\right)\!\!\bigg/\!\!\int_0^1 r^3 (\zeta-1)\dr,$$
for $G^\varepsilon$ and its inverse $(G^\varepsilon)^{-1}:(\eta,\hat{\omega},\hat{\phi},\zeta) \mapsto (\eta,\omega,\phi,\zeta)$.\qed

A simple calculation shows that the diffeomorphism $G$ transforms
$$u_z = g^\varepsilon(u)$$
into
\begin{equation}
u_z=\hat{g}^\varepsilon(u), \label{HE with linear BC}
\end{equation}
where $\hat{g}^\varepsilon: Y \rightarrow X$ is the smooth vector field defined by
$$
\hat{g}^\varepsilon(u)={{\mathrm d}G^\varepsilon}\,[(G^\varepsilon)^{-1}(u)]\,
(g^\varepsilon((G^\varepsilon)^{-1}(u))).
$$
Formula \eqref{HE with linear BC} represents Hamilton's equations for the Hamiltonian
system $(\hat{M},\Upsilon^\varepsilon,\hat{H}^\varepsilon)$, where
$$
\Upsilon^\varepsilon\big|_m(v_1,v_2)=
\Omega({{\mathrm d}G}^\varepsilon[(G^\varepsilon)^{-1}(m)]^{-1}(v_1),
{{\mathrm d}G}^\varepsilon[(G^\varepsilon)^{-1}(m)]^{-1}(v_2)), \quad m \in \hat{M}, \ v_1,v_2 \in T\hat{M}|_m,
$$
and
$$
\hat{H}^\varepsilon(n)=H^\varepsilon((G^\varepsilon)^{-1}(n)), \qquad n \in \hat{N}.
$$
The domain of the Hamiltonian vector field $v_{\hat{H}^\varepsilon}$ is
$$\DD(v_{\hat{H}^\varepsilon})=\{(\eta,\hat{\omega},\hat{\phi},\zeta) \in \hat{N}: r \hat{\phi}_r|_{r=1}=0\}$$
and
$v_{\hat{H}^\varepsilon}|_{n}=\hat{g}^\varepsilon(n)$ for any $n\in \DD(v_{\hat{H}^\varepsilon})$.

The next step is to verify that \eqref{HE with linear BC} satisfies the hypotheses
of Theorem \ref{Mielke's theorem} (with $\XX=X$), so that we obtain
a finite-dimensional reduced Hamiltonian system
$(\tilde{M}^\varepsilon,\tilde{\Gamma}^\varepsilon,\tilde{H}^\varepsilon)$.
We write \eqref{HE with linear BC} as
$$
u_z = L u + N^\varepsilon(u),
$$
in which $L={\mathrm d} v_{\hat{H}^0}[0]$ and verify the spectral hypotheses on $L$
by considering the operator $K: \DD(K) \subseteq X \rightarrow X$, where
\begin{equation}
K\begin{pmatrix} \eta \\ \omega \\ \phi \\ \zeta \end{pmatrix}
=
\begin{pmatrix}
\displaystyle\frac{1}{\beta_0}\left(\omega-\int_0^1 r^2\phi_r \dr\right)
\\[0.5em]
\displaystyle -2\int_0^1 r \zeta \dr - 2 \eta +(\alpha_0-\beta_0)\eta \\[0.5em]
\zeta + 2\eta \\[0.5em]
\displaystyle -\frac{1}{r}(r\phi_r)_r-\frac{2}{\beta_0}\left(\omega-\int_0^1 r^2\phi_r \dr\right)
\end{pmatrix}
\label{Formula for K}
\end{equation}
and
$$\DD(K)=\left\{(\eta,\omega,\phi,\zeta) \in Y: -r\phi_r -\frac{r^2}{\beta_0} \left(\omega-\int_0^1 r^2\phi_r\dr\right)\Bigg|_{r=1}\right\}$$
(the formal linearisation of $v_{H^0}$ at the origin); the formula $K=\mathrm{d}G^0[0]^{-1}L\mathrm{d}G^0[0]$ shows that the
spectral properties of $K$ and $L$ are identical.
It follows from Lemma \ref{lin} below that $L$ satisfies hypotheses (H1) and (H2);
hypothesis (H3) is clearly satisfied for an arbitrary value of $k$. Part (i) of Lemma \ref{lin}
is proved using the elementary theory of ordinary differential equations, while part
(ii) is established using arguments similar to those employed for other problems
treated using centre-manifold reduction (e.g.\ see Buffoni, Groves \& Toland
\cite[Proposition 3.2]{BuffoniGrovesToland96} or Groves \& Wahl\'{e}n \cite[Lemma 3.4]{GrovesWahlen07}).

\begin{lemma}\label{lin}
\hfill
\begin{list}{(\roman{count})}{\usecounter{count}}
\item
The spectrum $\sigma(L)$ of $L$ consists entirely of isolated
eigenvalues of finite algebraic multiplicity.
A complex number $\lambda$ is an eigenvalue of $L$ if and only if
$$\lambda J_0(\lambda)=(\gamma_0-\beta_0\lambda^2)J_1(\lambda),$$
where $\gamma_0=\alpha_0-\beta_0$.
(In particular, $0$ is an eigenvalue of $L$ and $\sigma(L) \cap \i {\mathbb R}$ is a finite set.)
\item
There exist real constants $C$, $s_0>0$ such that
$$
\|(L-\i sI)^{-1}\|_{\LL(X,X)}\leq\frac{C}{|s|}
$$
for each real number $s$ with $|s|>s_0$.
\end{list}
\end{lemma}

\begin{figure}
\centering
\includegraphics[scale=0.48]{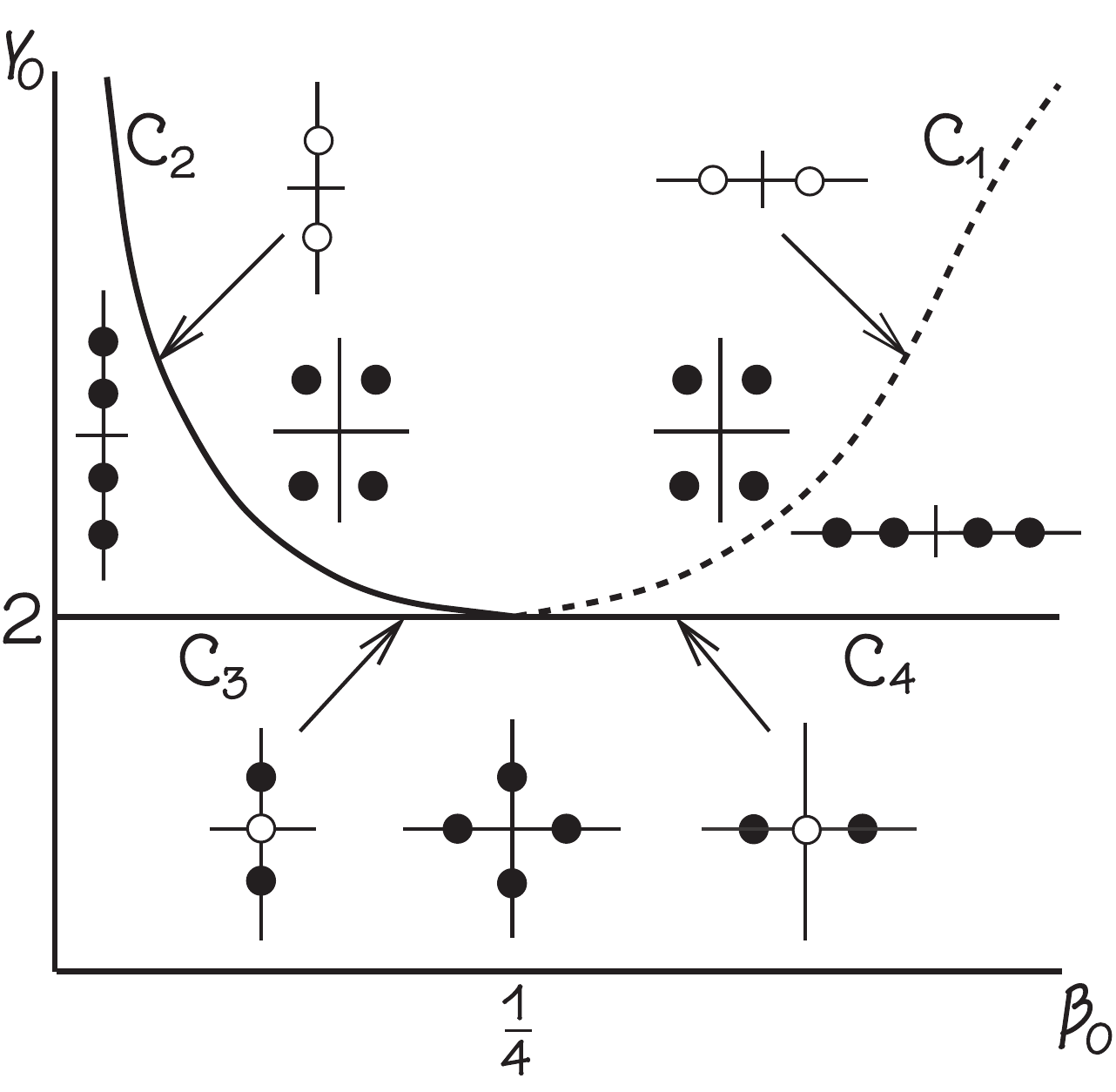}
{\it\caption{
Eigenvalues of $L$; solid and hollow dots denote respectively algebraically
simple and multiple eigenvalues. The curves $C_j$, $j=1,\ldots, 4$ consist of points in
$(\beta_0, \gamma_0)$ parameter space at which the qualitative nature of the eigenvalue picture changes.
\label{Bifurcation curves}}}
\end{figure}

According to Lemma \ref{lin}(i), a purely imaginary number $\lambda=\mathrm{i}s$ is an eigenvalue of $L$ if and only if
$$s I_0(s)=(\gamma_0+\beta_0 s^2)I_1(s).$$
Straightforward computations show that there are three critical curves
\[C_2 = \left\{(\beta_0,\gamma_0)=\left(\frac{1}{2}\left(1-\frac{I_0(s)I_2(s)}{I_1(s)^2}\right),\frac{1}{2}s^2\left(-1+\frac{I_0(s)^2}{I_1(s)^2}\right)\right)
:s \in (0,\infty) \right\}\]
and
\[C_3=\{(\beta_0,\gamma_0): \beta_0 < \tfrac{1}{4}, \gamma_0=2\}, \qquad
C_4=\{(\beta_0,\gamma_0): \beta_0 > \tfrac{1}{4}, \gamma_0=2\}\]
in the $(\beta_0,\gamma_0)$ parameter plane at which purely imaginary
eigenvalues of $L$ collide, together with a fourth curve
\[C_1=\left\{(\beta_0,\gamma_0)=\left(
\frac{1}{2}\left(1-\frac{J_0(k)J_2(k)}{J_1(k)^2}\right),\frac{k^2(J_0(k)^2+J_1(k)^2)}{2J_1(k)^2}\right): k \in (0,j_{1,1})\right\}\]
at which real eigenvalues collide (see Figure \ref{Bifurcation curves}). Here $J_0$, $J_1$, \ldots and $I_0$, $I_1$, \ldots
denote respectively the Bessel functions and modified Bessel functions of the first kind, and $j_{1,1}>0$ is the smallest zero of $J_1$.
Furthermore, $L$ has a geometrically simple zero
eigenvalue whose algebraic multiplicity is two for $\gamma_0 \neq 2$, four for $\gamma_0=2$, $\beta_0 \neq \frac{1}{4}$ and
and six for $(\beta_0,\gamma_0)=(\frac{1}{4},2)$.

The centre manifold  $\tilde{M}^\varepsilon$
is equipped with the single coordinate chart $\tilde{U}_1 \subset \XX_1$ and  coordinate
map $\pi: \tilde{M}^\varepsilon \rightarrow \tilde{U}_1$ defined by
$\pi^{-1}(u_1) = u_1 + r(u_1;\varepsilon)$. It is however more convenient to use
an alternative coordinate map for calculations.
We define the function
$\tilde{r}: \tilde{W}_1 \times \tilde{\Lambda} \rightarrow \tilde{U}_1 \times \tilde{U}_2$
with $\tilde{W}_1=\PP (G^\varepsilon)^{-1}(\tilde{U}_1\times \tilde{U}_2)$
(which in general has components in $\XX_1$ and $\XX_2$) by the formula
\begin{equation}
w_1 + \tilde{r}(w_1;\varepsilon) = (G^\varepsilon)^{-1}\big(w_1+r(w_1;\varepsilon)\big), \label{Definition of tilde r}
\end{equation}
where $\tilde{r}(0;0)=0$, $\mathrm{d}_1\tilde{r}[0;0]=0$, and equip $\tilde{M}^\varepsilon$
with the coordinate map $\hat{\pi}: \tilde{M}^\varepsilon \rightarrow \tilde{W}_1$ given by
$\hat{\pi}^{-1}(w_1) = w_1 + \tilde{r}(w_1;\varepsilon)$, so that
\begin{align}
\tilde{H}^\varepsilon(w_1) & = H^\varepsilon(w_1+\tilde{r}(w_1;\varepsilon)), \nonumber \\
\tilde{\Omega}^\varepsilon|_{w_1}(v_1,v_2) & =
\Omega(v_1+\mathrm{d}_1\tilde{r}[w_1;\varepsilon](v_1),
v_2+\mathrm{d}_1\tilde{r}[w_1;\varepsilon](v_2)) \nonumber  \\
& =\Omega(v_1,v_2) + O(|(\varepsilon,w_1)|) \label{Reduced 2 form}
\end{align}
as $(\varepsilon,w_1) \rightarrow 0$. Furthermore, using a parameter-dependent
version of Darboux's theorem (e.g.\ see
Buffoni \& Groves \cite[Theorem 4]{BuffoniGroves99}), we may assume that the remainder term
in \eqref{Reduced 2 form} vanishes identically.\pagebreak

We proceed by choosing a symplectic basis $\{f_0^1,\ldots,f_n^1,f_0^2,\ldots,f_n^2\}$ for the centre subspace of $K$ (so that $\Omega(f_i^1,f_i^2)=1$ for
$i=0,\ldots,n$ and the symplectic product of any other combination of these vectors is zero); here either $f_0^1$ or $f_0^2$ is the eigenvector
$(0,0,1,0)^\mathrm{T}$ corresponding to the zero eigenvalue of $K$. Using coordinates $q_0$, \ldots, $q_n$, $p_0$, \ldots, $p_n$, where
$$w_1 = q_0f_0^1 + q_1f_1^1 + \cdots + q_nf_n^1 +p_0f_0^2+p_1f_1^2 + \cdots +p_nf_n^2,$$
we find that $\tilde{\Omega}^\varepsilon$ is the canonical $2$-form.
Note that equations \eqref{HE 1}--\eqref{BC 2} are invariant under the transformation $\phi \mapsto \phi + c$, $c \in {\mathbb R}$,
and the quantity $\int_0^1 r \zeta \dr$ is conserved. This symmetry is inherited by the reduced system: one of the variables $q_0$, $p_0$ is cyclic
(that is, $\tilde{r}$ and $\tilde{H}^\varepsilon$ do not depend upon it), so that the other is conserved.

According to the
classical theory, the next step is to  lower the dimension of the reduced system by two
by setting the conserved variable to zero, solving the resulting decoupled system for $q_1$, \ldots, $q_n$, $p_1$, \ldots $p_n$,
and recovering the cyclic variable by quadrature; the lower-order system is typically studied using a canonical change
of variables which simplifies its Hamiltonian $\tilde{H}^\varepsilon|_{q_0=0}$
(a `normal-form' transformation). For our purposes it is convenient to use a normal-form
transformation \emph{before} lowering the order of the system since it can be `absorbed' into $\tilde{r}$ in the same way as the Darboux
transformation; this procedure greatly simplifies our later calculations. The following general result (whose proof is
based upon the method given by Bridges \& Mielke \cite[Theorem 4.3]{BridgesMielke95}) shows that this procedure is possible;
we assume for definiteness that $p_0$ is cyclic and use the construction by Elphick \cite{Elphick88} as our `usual' normal form.
The result is applied to the specific parameter regimes shown in Figure \ref{summary}(a) in Section \ref{Analysis of reduced systems}
below, where we denote the nonlinear part of the reduced Hamiltonian
vector field $v_{\tilde{H}^\varepsilon}$ by $P^\varepsilon(w_1)$.

\begin{theorem} \label{BM theorem}
Consider the $(n+1)$-degree-of-freedom Hamiltonian system
\begin{align}
& \dot{q}_i = \frac{\partial \tilde{H}^\varepsilon}{\partial p_i}, \quad \dot{p}_i = -\frac{\partial \tilde{H}^\varepsilon}{\partial q_i}, \qquad i=1,\ldots,n, \label{General NF 1} \\
& \dot{q}_0 = \frac{\partial \tilde{H}^\varepsilon}{\partial p_0}, \quad \dot{p}_0 = -\frac{\partial \tilde{H}^\varepsilon}{\partial q_0}, \label{General NF 2}
\end{align}
where $\tilde{H}^\varepsilon(q,p,q_0)=O(|(\varepsilon,q_0,q,p)||(q_0,q,p)|)$ and $p_0$ is cyclic (so that $q_0$ is conserved).

There exists a near-identity canonical change of variables $(q,p,q_0,p_0) \mapsto (Q,P,Q_0,P_0)$ with the properties that
$P_0$ is cyclic, $Q_0=q_0$ and the lower-order Hamiltonian system
$$
\dot{Q}_i = \frac{\partial \tilde{H}^\varepsilon}{\partial P_i}(Q,P,0), \quad \dot{P}_i = -\frac{\partial \tilde{H}^\varepsilon}{\partial Q_i}(Q,P,0), \qquad i=1,\ldots,n,$$
adopts its usual normal form. (Here, with a slight abuse of notation, we denote the transformed Hamiltonian by $\tilde{H}^\varepsilon(Q,P,Q_0)$.)
\end{theorem}

{\bf Proof.} Consider the $n$-degree of freedom Hamiltonian system
\begin{equation}
\dot{q}_i = \frac{\partial \tilde{H}^\varepsilon}{\partial p_i}, \quad \dot{p}_i = -\frac{\partial \tilde{H}^\varepsilon}{\partial q_i}, \qquad i=1,\ldots,n,
\label{Parameter-dependent system}
\end{equation}
in which $q_0$ and $\varepsilon$ are parameters. The standard theory asserts the existence of a canonical change of variables
\begin{align*}
Q & = q+ h_1^\varepsilon(q,p,q_0), \\
P & = p + h_2^\varepsilon(q,p,q_0)
\end{align*}
with
$$h_j^\varepsilon(q,p,q_0)=O(|(\varepsilon,q_0,q,p)||(q,p)|), \qquad j=1,2,$$
which converts \eqref{Parameter-dependent system} into its parameter-dependent normal form; note that
$$M_1^\mathrm{T}J_1M_1=J_1,$$
where
$$
M_1 = \begin{pmatrix} Q_q & Q_p \\ P_q & P_q \end{pmatrix} = \begin{pmatrix} I + \partial_q h_1^\varepsilon & \partial_p h_1^\varepsilon \\
\partial_q h_2^\varepsilon & I+\partial_p h_2^\varepsilon \end{pmatrix}, \qquad
J_1 = \begin{pmatrix} 0 & I \\ -I & 0 \end{pmatrix},
$$
and this condition may also be written as
\begin{equation}
(I+\partial_q h_1^\varepsilon)(I+\partial_p h_2^\varepsilon) - \partial q h_2^\varepsilon \partial_p h_1^\varepsilon = I.
\label{Condition 0}
\end{equation}

We seek a change of variable for \eqref{General NF 1}, \eqref{General NF 2} of the form
\begin{align*}
Q & = q + h_1^\varepsilon(q,p,q_0), \\
P & = p + h_2^\varepsilon(q,p,q_0), \\
Q_0 & = q_0, \\
P_0 & = p_0 + h_4^\varepsilon(q,p,q_0);
\end{align*}
the new function
$$h_4^\varepsilon(q,p,q_0)=O(|(\varepsilon,q_0,q,p)||(q,p)|)$$
is subject to the requirement that
$$M_2^\mathrm{T}J_2M_2 = J_2,$$
where
\begin{align*}
M_2 & = \begin{pmatrix} Q_q & Q_p & Q_{q_0} & Q_{p_0} \\
P_q & P_p & P_{q_0} & P_{p_0} \\
Q_{0q} & Q_{0p} & Q_{0q_0} & Q_{0p_0} \\
P_{0q} & P_{0p} & P_{0q_0} & P_{0p_0}\end{pmatrix}
= \begin{pmatrix} I + \partial_q h_1^\varepsilon & \partial_p h_1^\varepsilon & \partial_{q_0} h_1^\varepsilon & 0 \\
\partial_q h_2^\varepsilon & I + \partial_p h_2^\varepsilon & \partial_{q_0}h_2^\varepsilon & 0 \\
0 & 0 & 1 & 0 \\
\partial_q h_4^\varepsilon & I + \partial_p h_4^\varepsilon & \partial_{q_0}h_4^\varepsilon & 1 \end{pmatrix}, \\
J_2 & = \begin{pmatrix} 0 & I & 0 & 0 \\
-I & 0 & 0 & 0 \\
0 & 0 & 0 & 1 \\
0& 0 & -1 & 0
\end{pmatrix},
\end{align*}
and this condition may be written as
\begin{align}
(I+\partial_q h_1^\varepsilon)(I+\partial_ph_2^\varepsilon)- \partial_q h_2^\varepsilon \partial_p h_1^\varepsilon & = I, \label{Condition 1} \\
(I+\partial_q h_1^\varepsilon)\partial_{q_0}h_2^\varepsilon - \partial_q h_2^\varepsilon \partial_{q_0}h_1^\varepsilon & = \partial_q h_4^\varepsilon, \label{Condition 2} \\
\partial_p h_1^\varepsilon \partial_{q_0}h_2^\varepsilon - (I+\partial_p h_2^\varepsilon)\partial_{q_0}h_1^\varepsilon & = \partial_p h_4^\varepsilon. \label{Condition 3}
\end{align}
It is possible to find $h_4^\varepsilon$ satisfying these conditions since the compatibility condition for 
\eqref{Condition 2}, \eqref{Condition 3} is the derivative of \eqref{Condition 1} with respect to $q_0$, and \eqref{Condition 1} is automatically
satisfied because of \eqref{Condition 0}.\qed

\section{The reduced Hamiltonian systems} \label{Analysis of reduced systems}

\subsection{Homoclinic bifurcation at $C_4$} \label{Region I}

At each point of the curve $C_4$ in Figure \ref{Bifurcation curves} two real eigenvalues become purely imaginary by
colliding at the origin and increasing the algebraic multiplicity of the zero eigenvalue
from two to four. This resonance
is associated with the bifurcation of a branch of homoclinic solutions into the region
with real eigenvalues (the parameter regime marked I in Figure \ref{summary}.
Let us therefore fix reference values $(\beta_0,\gamma_0) \in C_4$, so that
$\beta_0>\frac{1}{4}$, $\alpha_0=2+\beta_0$, and introduce a bifurcation parameter
by choosing $(\varepsilon_1,\varepsilon_2) = (0,\mu)$, where $0 < \mu \ll 1$.

The four-dimensional centre subspace of $K$ is spanned by the generalised eigenvectors
$$
e_1 = \begin{pmatrix} 0 \\ 0 \\ 1 \\ 0 \end{pmatrix}, \quad
e_2 =\begin{pmatrix} \frac{1}{2} \\ 0 \\ 0 \\ 0 \end{pmatrix}, \quad
e_3 = \begin{pmatrix} 0 \\ \frac{1}{2}(\beta_0-\frac{1}{4}) \\ - \frac{1}{4}r^2 + A_4 \\ 0 \end{pmatrix}, \quad
e_4 = \begin{pmatrix} \frac{1}{4}(\beta_0-\frac{1}{2}) + \frac{1}{2}A_4 \\ 0 \\ 0 \\ - \frac{1}{4}r^2 - \frac{1}{2}(\beta_0-\frac{1}{2}) \end{pmatrix},
$$
where
$A_4=-(\beta_0-\tfrac{1}{4})^{-1}\left(\tfrac{1}{24}+\tfrac{1}{4}\beta_0(\beta_0-1)\right)$ has been chosen so that
$Ke_1=0$, $Ke_j=e_{j-1}$ for $j=2,3,4$,
$$\Omega(e_1,e_4)=-\tfrac{1}{4}\left(\beta_0-\tfrac{1}{4}\right), \qquad
\Omega(e_2,e_3)=\tfrac{1}{4}\left(\beta_0-\tfrac{1}{4}\right)$$
and the symplectic product of any other combination of the vectors $e_1$, \ldots $e_4$ is zero.
Writing
$$w_1=q_0 f_4 + p_0 f_1 +q f_2 + p f_3, \qquad f_i:=\tfrac{1}{2}\left(\beta_0-\tfrac{1}{4}\right)^{-1/2}e_i,$$
we therefore find that $q_0$, $q$, $p_0$ and $p$ are canonical coordinates for the reduced Hamiltonian system,
which has the cyclic variable $p_0$ and reverser $S:(q_0,q,p_0,p) \mapsto (q_0,q,-p_0,-p)$; with a slight abuse
of notation we abbreviate $\tilde{H}^{\varepsilon}|_{(\varepsilon_1,\varepsilon_2)=(0,\mu)}$ to $\tilde{H}^\mu$.

The usual normal-form theory for the two-dimensional system with Hamiltonian $\tilde{H}^\mu(q,p,0)$ asserts that,
after a canonical change of variables,
$$\tilde{H}^\mu(q,p,0) = \tfrac{1}{2}p^2 + \tilde{H}_\mathrm{NF}^0(q,\mu)+ O(|(q,p)|^2|(\mu,q,p)|^{n_0}),$$
where
$\tilde{H}_\mathrm{NF}^0(q,\mu)$ is a polynomial of order $n_0+1$ in $(q,\mu)$
with
\[\tilde{H}_\mathrm{NF}^0(q,\mu)=O(|q|^2|(\mu,q)|).\]
It follows that, after a canonical change of variables,
$$\tilde{H}^\mu(q,p,q_0) = \tfrac{1}{2}p^2 -q q_0+\tilde{H}^\mu_\mathrm{nl}(q,p,q_0)$$
with
$$
\tilde{H}^\mu_\mathrm{nl}(q,p,q_0)= \tilde{H}_\mathrm{NF}(q,q_0,\mu)+ \tilde{H}_\mathrm{r}(q,q_0,\mu)+ O(|(q,p,q_0)|^2|(\mu,q,p,q_0)|^{n_0});$$
here $\tilde{H}_\mathrm{NF}(q,q_0,\mu)$ is a polynomial of order $n_0+1$
with
\[\tilde{H}_\mathrm{NF}(q,q_0,\mu)=O(|q|^2|(\mu,q,q_0)|)\]
and $\tilde{H}_\mathrm{NF}(q,0,\mu)=\tilde{H}_\mathrm{NF}^0(q,\mu)$,
and $\tilde{H}_\mathrm{r}(q,q_0,\mu)$ is an affine function of its first argument which satisfies
$$\tilde{H}_\mathrm{r}(q,q_0,\mu)=O(|(q,q_0)||q_0||(\mu,q,q_0)|).$$
Note that
$$P^\mu(q,p,q_0)\!=\! -\partial_q \tilde{H}_\mathrm{nl}^\mu(q,p,q_0)f_3- \partial_{q_0} \tilde{H}_\mathrm{nl}^\mu(q,p,q_0)f_1.
$$

Writing
\begin{align*}
\tilde{H}_3^0(q,p,q_0) & = c_1 q^3 + c_2 q^2 q_0 + c_3 qq_0^2 + c_4 q_0^3, \\
\tilde{H}_2^1(q,p,q_0) & = c_1^1 q^2 + c_2^1 qq_0 + c_3^1 q_0^2,
\end{align*}
where $\mu^j\tilde{H}^j_k(q,p,q_0)$ denotes
the part of the Taylor expansion of $\tilde{H}^\mu(q,p,q_0)$
which is homogeneous of order $j$ in $\mu$ and $k$ in $(q,p,q_0)$,
one finds that
$$c_1 = \tfrac{1}{6}\left(\beta_0-\tfrac{1}{4}\right)^{-3/2}(\alpha_0 m_1^\prime(1)-6),
\qquad
c_1^1 =-\tfrac{1}{2}\left(\beta_0-\tfrac{1}{4}\right)^{-1}$$
(see Appendix (i)).
Setting $q_0=0$ and introducing scaled variables
$$Z=\mu^{1/2}(\beta_0-\tfrac{1}{4})^{-1/2}z, \quad q(z)=\mu (\beta_0-\tfrac{1}{4})^{1/2}Q(Z),\quad p(z)=\mu^{3/2}P(Z),$$
yields
$$
\tilde{H}^\mu(q,p,0)=\mu^3 \left[\tfrac{1}{2}P^2 - \tfrac{1}{2}Q^2 + \tfrac{1}{3}\check{c}_1 Q^3\right] + O(\mu^{7/2}),
$$
where
$$\check{c}_1= \tfrac{1}{2}(\alpha_0 m_1^\prime(1)-6),$$
and the lower-order Hamiltonian system
\begin{align}
\dot{Q} & = P+O(\mu^{1/2}), \label{Scaled reduced system 1}\\
\dot{P} & = Q - \check{c}_1 Q^2 + O(\mu^{1/2}), \label{Scaled reduced system 2}
\end{align}
which is reversible with reverser $S:(Q,P) \mapsto (Q,-P)$. Suppose $\check{c}_1 \neq 0$.
In the limit $\mu=0$ equations \eqref{Scaled reduced system 1}, \eqref{Scaled reduced system 2}
are equivalent to the single equation
$$\partial_Z^2u-u+u^2=0$$
for the variable $u=\check{c}_1Q$.

Let us now suppose that $m_1^\prime(1)$ is close to the critical value $6\alpha_0^{-1}$ and introduce a second bifurcation
parameter $\kappa$ by setting
$$m_1^\prime(1)=\alpha_0^{-1}(6+\kappa)$$
and observing that
\begin{align*}
\tilde{r}(q,p,q_0;\mu,\kappa) &= O(|(q,p,q_0)||(\mu,q,p,q_0)|) + O(|\kappa||(q,p,q_0)|^2), \\
\tilde{H}^{\mu,\kappa}(q,p,q_0) &= O(|(q,p,q_0)|^2|(\mu,q,p,q_0)|) + O(|\kappa||(q,p,q_0)|^3)
\end{align*}
(with a slight change of notation).
Writing
$$
\tilde{H}_4^{0,0}(q,p,q_0) = d_1 q^4 + d_2 q^3 q_0 + d_3 q^2q_0^2 + d_4 qq_0^3 + d_5q_0^4,
$$
where
$\mu^i\kappa^j\tilde{H}^{i,j}_k(q,p,q_0)$ denotes
the part of the Taylor expansion of $\tilde{H}^{\mu,\kappa}(q,p,q_0)$
which is homogeneous of order $i$ in $\mu$, $j$ in $\kappa$ and $k$ in $(q,p,q_0)$,
one finds that
$$d_1 = \tfrac{1}{24}(\beta_0-\tfrac{1}{4})^{-2}(12-\alpha_0 m_1^{\prime\prime}(1))$$
(see Appendix (ii)). Setting $q_0=0$, introducing scaled variables
$$Z=\mu^{1/2}(\beta_0-\tfrac{1}{4})^{-1/2}z, \quad q(z)=\mu^{1/2} (\beta_0-\tfrac{1}{4})^{1/2}Q(Z),\quad p(z)=\mu P(Z)$$
and writing
$$\kappa= 2\mu^{1/2}\check{\kappa},$$
thus yields
$$
\tilde{H}^{\mu,\kappa}(q,p,0)=\mu^2 \left[\tfrac{1}{2}P^2 - \tfrac{1}{2}Q^2 + \tfrac{1}{3}\check{\kappa}Q^3 + \tfrac{1}{4}\check{d}_1 Q^4\right] + O(\mu^{5/2}),
$$
where
$$\check{d}_1=\tfrac{1}{6}(12-\alpha_0 m_1^{\prime\prime}(1)),$$
and the lower-order Hamiltonian system
\begin{align}
\dot{Q} & = P+O(\mu^{1/2}), \label{Scaled reduced system 3}\\
\dot{P} & = Q - \check{\kappa} Q^2 -\check{d}_1 Q^3+O(\mu^{1/2}), \label{Scaled reduced system 4}
\end{align}
which is of course reversible with reverser $S:(Q,P) \mapsto (Q,-P)$. Suppose that $\check{d}_1>0$. In the limit $(\mu,\check{\kappa})=0$
equations \eqref{Scaled reduced system 3}, \eqref{Scaled reduced system 4}
are equivalent to the single equation
$$\partial_Z^2u-u+u^3=0$$
for the variable $u=\check{d}_1^{1/2}Q$.

\begin{figure}[h]
\centering
\includegraphics[scale=0.8]{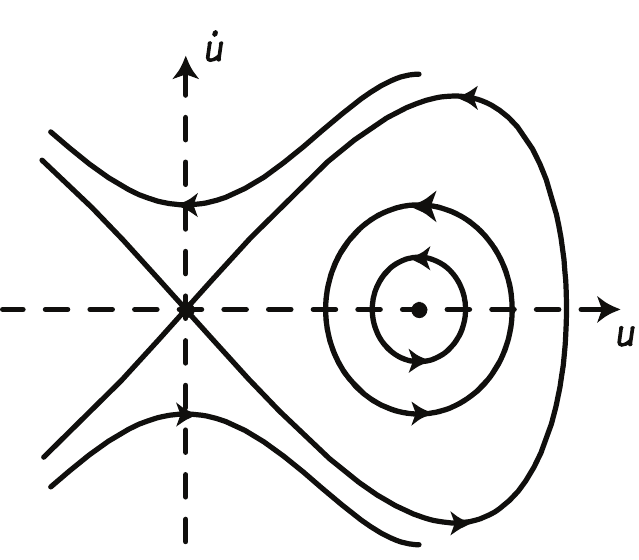}\hspace{2cm}
\includegraphics[scale=0.8]{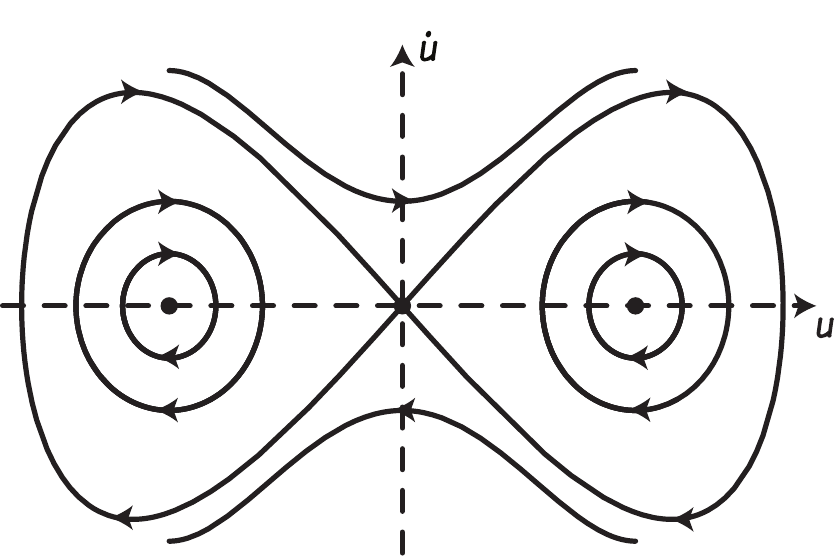}
{\it
\caption{Phase portrait of equation \eqref{KdV-style eqn} for even (left) and odd (right) values of $m$} \label{C1 phase portrait}}
\vspace*{8mm}
\includegraphics[scale=0.7]{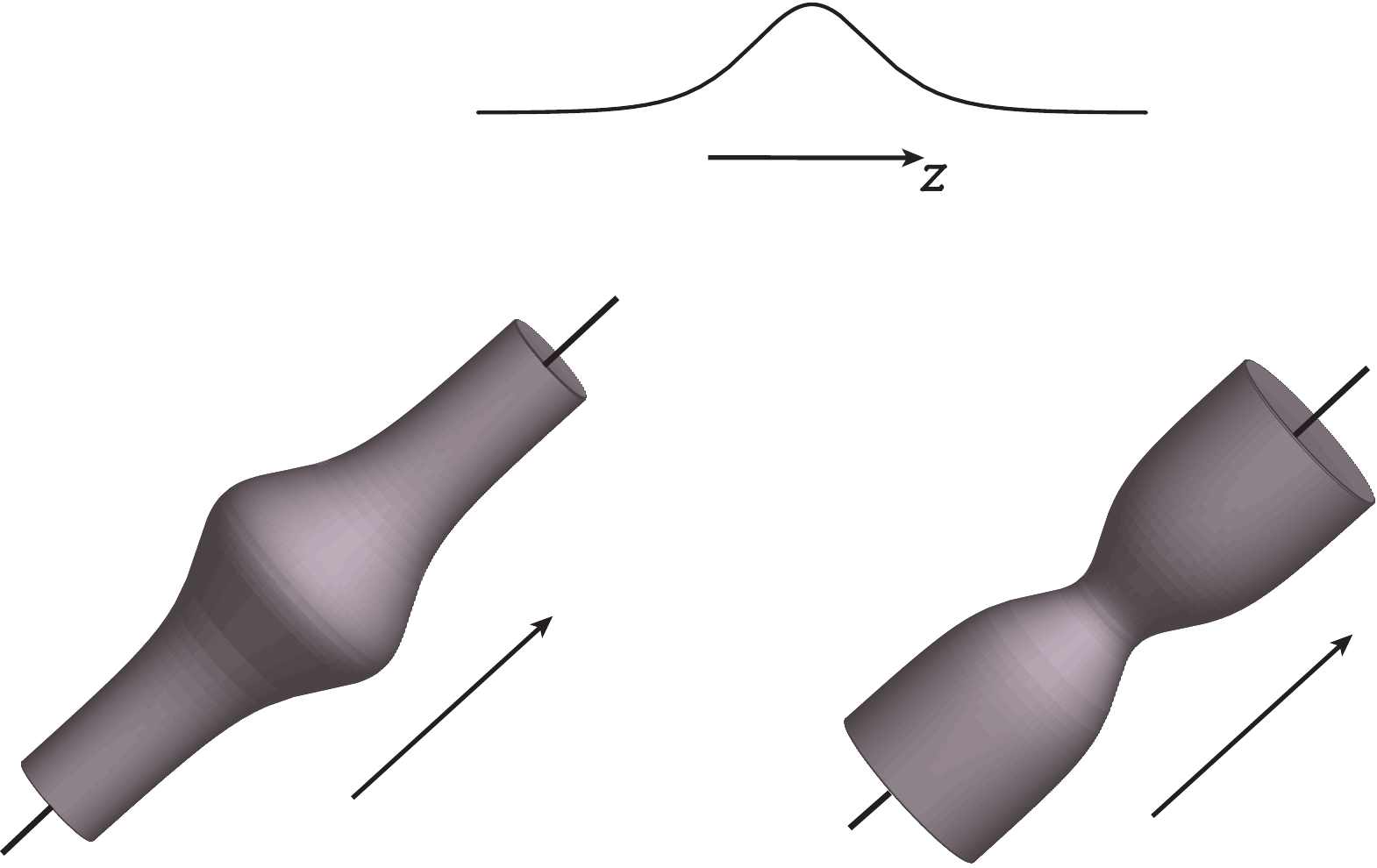}
{\it
\caption{A solitary wave of elevation (left) and depression (right) generated by a
homoclinic solution (top) in region I} \label{C4 waves}}
\end{figure}

The phase portrait of the equation
\begin{equation}
\ddot{u}-u+u^m=0 \label{KdV-style eqn}
\end{equation}
for a fixed natural number $m$ (which is a travelling-wave version of the generalised Korteweg-de Vries equation)
is readily obtained by elementary calculations
and is sketched in Figure \ref{C1 phase portrait}; the homoclinic orbits are of particular interest.

\begin{lemma} \label{KdV homoclinics} $ $
\begin{itemize}
\item[(i)]
Suppose that $m$ is even. Equation \eqref{KdV-style eqn} has precisely one homoclinic solution $h$ (up to translations). This
solution is positive and symmetric, and monotone increasing to the left, monotone decreasing to the right of its point of symmetry.
\item[(ii)]
Suppose that $m$ is odd. Equation \eqref{KdV-style eqn} has precisely two homoclinic solutions
$\pm h$, where $h$ is symmetric, and monotone increasing to the left, monotone decreasing to the right of its point of symmetry.
\end{itemize}
In both cases the homoclinic solutions intersect the symmetric section $\{\dot{u}=0\}$ in the two-dimensional
phase space $\{(u,\dot{u}) \in {\mathbb R}^2\}$ transversally.
\end{lemma}

A familiar argument shows that Lemma \ref{KdV homoclinics}(i) also applies to \eqref{Scaled reduced system 1}, \eqref{Scaled reduced system 2}
for small, positive values of $\mu$, while Lemma \ref{KdV homoclinics}(ii) applies to \eqref{Scaled reduced system 3}, \eqref{Scaled reduced system 4}
for small, positive values of $\mu$ and small, values of $\check{\kappa}$ (that is, small values of $\kappa$); the qualitative statements
apply to the variable $\check{c}_1Q$ or $\check{d}_1^{1/2}Q$. The homoclinic orbits at $\mu=0$
(and $\check{\kappa}=0$) intersect the symmetric section $\mbox{Fix}\, R=\{P=0\}$ transversally, and these orbits therefore persist (as small, uniform perturbations
of their limits) for small, positive values of $\mu$ (and small values of $\check{\kappa}$).

Altogether we have established the existence of a symmetric, monotonically decaying solitary wave of
depression for $m_1^\prime(1) < 6\alpha_0^{-1}$ and elevation for $m_1^\prime(1) > 6\alpha_0^{-1}$;
the corresponding ferrofluid surface $\{r=1+\eta(z)\}$ is obtained from the homoclinic
solution of \eqref{Scaled reduced system 1}, \eqref{Scaled reduced system 2}
by the formula
$$\eta(z)=\tfrac{1}{2}\mu(\beta_0-\tfrac{1}{4})^{1/2}Q\left(\mu^{1/2}(\beta_0-\tfrac{1}{4})^{-1/2}z\right) + O(\mu^{3/2}).$$
Furthermore, a pair of symmetric, monotonically decaying solitary waves exists for small values of 
$m_1^\prime(1)-6\alpha_0^{-1}$ provided that $m_1^{\prime\prime}(1) <12\alpha_0^{-1}$;
one is a wave of depression, the other a wave of elevation. The corresponding ferrofluid surface $\{r=1+\eta(z)\}$ is obtained from a homoclinic
solution of \eqref{Scaled reduced system 3}, \eqref{Scaled reduced system 4} by the formula
$$\eta(z)=\tfrac{1}{2}\mu^{1/2}(\beta_0-\tfrac{1}{4})^{1/2}Q\left(\mu^{1/2}(\beta_0-\tfrac{1}{4})^{-1/2}z\right) + O(\mu^{3/2}).$$
(A more detailed analysis of a codimension-two bifurcation of this kind was given by Kirrmann \cite[\S 4]{Kirrmann91} in the context
of two-layer fluid flow.)
Figure \ref{C4 waves} shows a sketch of the ferrofluid surface corresponding to solitary waves of the present type.

\subsection{Homoclinic bifurcation at $C_1$} \label{Region II}

At each point of the curve $C_1$ in Figure \ref{Bifurcation curves} two pairs of real eigenvalues become
complex by colliding at non-zero points on the real axis. Of particular interest here is the local part of $C_1$
near the point $(\beta_0,\gamma_0)=(\frac{1}{4},2)$ (which is given by $\beta_0=\frac{1}{4}+\frac{1}{48}\mu^2 + O(\mu^4)$,
$\gamma_0=2+\frac{1}{96}\mu^4 + O(\mu^6)$ for $0 < \mu \ll 1$)
since we can access this curve using the
centre-manifold technique. To this end we choose $\beta_0=\frac{1}{4}$, $\alpha_0=\frac{9}{4}$ and
\begin{equation}
\varepsilon_1=\mu_1, \quad \varepsilon_2 = \mu_1+\mu_2, \qquad 
\mu_1 = \tfrac{1}{48}(1+\delta)\mu^2, \quad \mu_2 = \tfrac{1}{96}\mu^4.
\label{C1 parameter choice}
\end{equation}
Notice that $\mu$ indicates the distance in $(\beta_0,\gamma_0)$ parameter space from the point $(\frac{1}{4},2)$,
while $\delta$ plays the role of a bifurcation parameter (varying $\delta$ through zero
from above we cross the critical curve $C_1$ from above); the
parameter regime marked II in Figure \ref{summary} corresponds
to small, positive values of $\delta$ and $\mu$.

The six-dimensional centre subspace of $K$ is spanned by the generalised eigenvectors
$$
e_1 = \begin{pmatrix} 0 \\ 0 \\ 1 \\ 0 \end{pmatrix}, \quad
e_2 =\begin{pmatrix} \frac{1}{2} \\ 0 \\ 0 \\ 0 \end{pmatrix}, \quad
e_3 = \begin{pmatrix} 0 \\ 0 \\ \frac{3}{32} - \frac{1}{4}r^2 \\ 0 \end{pmatrix}, \quad
e_4 = \begin{pmatrix} -\frac{1}{64} \\ 0 \\ 0 \\ \frac{1}{8}-\frac{1}{4}r^2 \end{pmatrix},
$$
$$
e_5 = \begin{pmatrix} 0 \\ - \frac{1}{192} \\ \frac{87}{10240}-\frac{3}{128}r^2+\frac{1}{64}r^4 \\ 0 \end{pmatrix}, \quad
e_6 = \begin{pmatrix} -\frac{33}{20480} \\ 0 \\ 0 \\ \frac{3}{256}-\frac{3}{128}r^2+\frac{1}{64} r^4 \end{pmatrix}
$$
where
$Ke_1=0$, $Ke_j=e_{j-1}$ for $j=2,\ldots,6$,
$$\Omega(e_1,e_6)=\tfrac{1}{384}, \qquad \Omega(e_2,e_5)=-\tfrac{1}{384}, \qquad \Omega(e_3,e_4)=\tfrac{1}{384}$$
and the symplectic product of any other combination of the vectors $e_1$, \ldots $e_6$ is zero.
Writing
$$w_1=q_0f_1+p_0f_6+q_1f_5+p_1f_2+q_2f_3+p_2f_4, \qquad f_i:=8\sqrt{6}e_i,$$
we therefore find that $q_0$, $q_1$, $q_2$, $p_0$, $p_1$ and $p_2$ are canonical coordinates for the reduced Hamiltonian system,
which has the cyclic variable $q_0$ and reverser $S:(q_0,q_1,q_2,p_0,p_1,p_2) \mapsto$\linebreak$(-q_0,-q_1,-q_2,p_0,p_1,p_2)$;
with a slight abuse of notation we abbreviate $\tilde{H}^{\varepsilon}|_{(\varepsilon_1,\varepsilon_2)=(\mu_1,\mu_1+\mu_2)}$ to $\tilde{H}^{\mu_1,\mu_2}$.

The usual normal-form theory for the four-dimensional system with Hamiltonian $\tilde{H}^{\mu_1,\mu_2}(q,p,0)$,
where $q=(q_1,q_2)$, $p=(p_1,p_2)$ asserts that,
after a canonical change of variables,
$$
\tilde{H}^{\mu_1,\mu_2}(q,p,0) = \tfrac{1}{2}p_2^2 -q_1q_2 + \tilde{H}_\mathrm{NF}^0(q,p,\mu_1,\mu_2)
+ O(|(q,p)|^2|(\mu_1,\mu_2,q,p)|^{n_0}),
$$
where $\tilde{H}_\mathrm{NF}^0(q,p,\mu_1,\mu_2)$ is a polynomial of order $n_0+1$ which depends upon $q_1$, $q_2$, $p_1$, $p_2$ through the combinations
$$p_1, \quad q_2^2-2p_1p_2, \quad q_2^3+3p_1^2q_1-3p_1p_2q_1, \quad -8p_1p_2^3+3p_2^2q_2^2-9p_1^2q_1^2-6q_1q_2^3+18p_1p_2q_1q_2$$
and satisfies
\[\tilde{H}_\mathrm{NF}^0(q,p,\mu_1,\mu_2)=O(|(q,p)|^2|(\mu_1,\mu_2,q,p)|).\]
It follows that, after a canonical change of variables,
$$
\tilde{H}^{\mu_1,\mu_2}(q,p,p_0) = \tfrac{1}{2}p_2^2 - q_1q_2  + p_0p_1 + \tilde{H}^{\mu_1,\mu_2}_\mathrm{nl}(q,p,p_0)
$$
with
$$\tilde{H}^{\mu_1,\mu_2}_\mathrm{nl}(q,p,p_0)
=\tilde{H}_\mathrm{NF}(q,p,p_0,\mu_1,\mu_2)
+\tilde{H}_\mathrm{r}(q,p,p_0,\mu_1,\mu_2)
+ O(|(q,p,p_0)|^2|(\mu_1,\mu_2,q,p,p_0)|^{n_0});
$$
here $\tilde{H}_\mathrm{NF}(q,p,p_0,\mu_1,\mu_2)$ is a polynomial of order $n_0+1$
which depends upon $q_1$, $q_2$, $p_1$, $p_2$ through the above combinations and satisfies
$$\tilde{H}_\mathrm{NF}(q,p,q_0,\mu_1,\mu_2)\!=\!O(|(q,p)|^2|(\mu_1,\mu_2,q,p,p_0)|),$$
and $\tilde{H}_\mathrm{NF}(q,p,0,\mu_1,\mu_2)\!=\!\tilde{H}_\mathrm{NF}^0(q,p,\mu_1,\mu_2)$,
and $\tilde{H}_\mathrm{r}(q,p,p_0,\mu_1,\mu_2)$ is an affine function of its first two arguments which satisfies
$$\tilde{H}_\mathrm{r}(q,p,p_0,\mu_1,\mu_2)=O(|(q,p,p_0)||p_0||(\mu_1,\mu_2,q,p,p_0)|).$$
Note that
\begin{align*}
P^{\mu_1,\mu_2}(q,p,q_0) & = \partial_{p_1} \tilde{H}_\mathrm{nl}^{\mu_1,\mu_2}(q,p,p_0) f_5 + \partial_{p_2} \tilde{H}_\mathrm{nl}^{\mu_1,\mu_2}(q,p,p_0) f_3 \\
& \qquad\mbox{}-\partial_{q_1} \tilde{H}_\mathrm{nl}^{\mu_1,\mu_2}(q,p,p_0)f_2- \partial_{q_2} \tilde{H}_\mathrm{nl}^{\mu_1,\mu_2}(q,p,p_0)f_4
+\partial_{p_0} \tilde{H}_\mathrm{nl}^{\mu_1,\mu_2}(q,p,p_0)f_1.
\end{align*}

Writing
\begin{align*}
\tilde{H}_3^{0,0}(q,p,p_0) & = c_1 p_1^3 + c_2 p_0 p_1^2 + c_3 p_0^2 p_1 + c_4 p_0^3 + c_5p_1(q_2^2-2p_1p_2) + c_6p_0(q_2^2-2p_1p_2)+c_7p_0^2p_2, \\
\tilde{H}_2^{1,0}(q,p,p_0) & = c_1^{1,0}p_1^2 + c_2^{1,0}p_0p_1+c_3^{1,0}p_0^2+c_4^{1,0}(q_2^2-2p_1p_2)+c_5^{1,0}p_0p_2, \\
\tilde{H}_2^{2,0}(q,p,p_0) & = c_1^{2,0}p_1^2 + c_2^{2,0}p_0p_1+c_3^{2,0}p_0^2+c_4^{2,0}(q_2^2-2p_1p_2)+c_5^{2,0}p_0p_2, \\
\tilde{H}_2^{0,1}(q,p,p_0) & = c_1^{0,1}p_1^2 + c_2^{0,1}p_0p_1+c_3^{0,1}p_0^2+c_4^{0,1}(q_2^2-2p_1p_2)+c_5^{0,1}p_0p_2,
\end{align*}
where $\mu_1^i\mu_2^j\tilde{H}^{i,j}_k(q,p,p_0)$ denotes
the part of the Taylor expansion of $\tilde{H}^{\mu_1,\mu_2}(q,p,q_0)$
which is homogeneous of order $i$ in $\mu_1$, $j$ in $\mu_2$ and $k$ in $(q,p,p_0)$,
one finds that
$$c_1 = 48\sqrt{6}(3 m_1^\prime(1)-8),
\quad
c_1^{1,0}=0,
\quad
c_4^{1,0}=-16,
\quad
c_1^{2,0}=512,
\quad
c_1^{0,1} = -48
$$
(see Appendix (iii)).
Setting $p_0=0$, choosing $\mu_1$, $\mu_2$ according to \eqref{C1 parameter choice} and
introducing the scaled variables
$$Z=\mu z, 
\quad q_1(z)=\mu^7 Q_1(Z), \quad q_2(z)=\mu^5 Q_2(Z),
\quad p_1(z)=\mu^4 P_1(Z), \quad p_2(z)=\mu^6 P_2(Z),
$$
thus yields
\begin{eqnarray*}
\lefteqn{\tilde{H}^{\mu_1,\mu_2}(q,p,0)} \\
& &
=\mu^{12}\left[\tfrac{1}{2}P_2^2-\tfrac{1}{2}P_1^2-Q_1Q_2-\tfrac{1}{3}(1+\delta)(Q_2^2-2P_1P_2)+\tfrac{2}{9}(1+\delta)^2P_1^2 + c_1 P_1^3\right]
+O(\mu^{13})
\end{eqnarray*}
and the lower-order Hamiltonian system
\begin{align}
Q_{1Z} & = -P_1 + \tfrac{2}{3}(1+\delta)P_2 + \tfrac{4}{9}(1+\delta)^2P_1 + 3 c_1P_1^2+ O(\mu), \label{C1 system 1} \\
Q_{2Z} & = P_2  + \tfrac{2}{3}(1+\delta)P_1 + O(\mu), \label{C1 system 2} \\
P_{1Z} & = Q_2+ O(\mu), \label{C1 system 3} \\
P_{2Z} & = Q_1 + \tfrac{2}{3}(1+\delta)Q_2+ O(\mu), \label{C1 system 4}
\end{align}
which is reversible with reverser $S:(Q,P) \mapsto (-Q,P)$. Suppose $c_1 \neq 0$.
In the limit $\mu=0$
equations \eqref{C1 system 1}--\eqref{C1 system 4} are equivalent to the single
fourth-order ordinary differential equation
$$
\partial_Z^4 u-2(1+\delta)\partial_Z^2u+u -u^2 =0
$$
for the variable $u=3c_1P_1$.

Let us now suppose that $m_1^\prime(1)$ is close to the critical value $\frac{8}{3}$ and introduce a further bifurcation
parameter $\kappa$ by setting
$$m_1^\prime(1)=\tfrac{1}{3}(8+\kappa)$$
and observing that
\begin{align*}
\tilde{r}(q,p,p_0;\mu_1,\mu_2,\kappa) &= O(|(q,p,p_0)||(\mu_1,\mu_2,q,p,p_0)|) + O(|\kappa|(q,p,p_0)|^2), \\
\tilde{H}^{\mu_1,\mu_2,\kappa}(q,p,p_0) &= O(|(q,p,p_0)|^2|(\mu_1,\mu_2,q,p,p_0)|) + O(|\kappa||(q,p,p_0)|^3)
\end{align*}
(with a slight change of notation).
Writing
\begin{align*}
\tilde{H}_4^{0,0,0}(q,p,p_0) & = d_1p_1^4 + d_2 p_1^3p_0 + d_3 p_1^2p_0^2+ d_4 p_1p_0^3 + d_5 p_0^4
+d_6(q_2^2-2p_1p_2)^2 \\
& \qquad\mbox{} + d_7 p_1^2(q_2^2-2p_1p_2)
+d_8 p_0^2(q_2^2-2p_1p_2) + d_9p_0p_1 (q_2^2-2p_1p_2) \\
& \qquad\mbox{} +d_{10}(-8p_1p_2^3+3p_2^2q_2^2-9p_1^2q_1^2-6q_1q_2^3+18p_1p_2q_1q_2) + d_{11}p_2p_0^3
\end{align*}
where
$\mu_1^i\mu_2^j\kappa^k\tilde{H}^{i,j}_\ell(q,p,q_0)$ denotes
the part of the Taylor expansion of $\tilde{H}^{\mu_1,\mu_2}(q,p,q_0)$
which is homogeneous of order $i$ in $\mu_1$, $j$ in $\mu_2$, $k$ in $\kappa$ and $\ell$ in $(q,p,q_0)$,
one finds that
$$d_1 = 864\left(\sdfrac{1264}{75}-m_1^{\prime\prime}(1)\right)$$
(see Appendix (iv)).
Setting $p_0=0$, choosing $\mu_1$, $\mu_2$ according to \eqref{C1 parameter choice},
introducing the scaled variables
$$Z=\mu z, 
\quad q_1(z)=\mu^5 Q_1(Z), \quad q_2(z)=\mu^3 Q_2(Z),
\quad p_1(z)=\mu^2 P_1(Z), \quad p_2(z)=\mu^4 P_2(Z),
$$
and writing
$$\kappa = \sdfrac{1}{144\sqrt{6}}\check{\kappa}\mu^2$$
thus yields
\begin{eqnarray*}
\lefteqn{\tilde{H}^{\mu_1,\mu_2,\kappa}(q,p,0)} \\
& &\hspace{-5mm}
=\!\mu^8\!\left[\tfrac{1}{2}P_2^2-\tfrac{1}{2}P_1^2-Q_1Q_2-\tfrac{1}{3}(1+\delta)(Q_2^2-2P_1P_2)+\tfrac{2}{9}(1+\delta)^2P_1^2
+\tfrac{1}{3}\check{\kappa}P_1^3 + d_1 P_1^4\right]
\!+\!O(\mu^9)
\end{eqnarray*}
and the lower-order Hamiltonian system
\begin{align}
Q_{1Z} & = -P_1 + \tfrac{2}{3}(1+\delta)P_2 + \tfrac{4}{9}(1+\delta)^2P_1 + \check{\kappa} P_1^2+4d_1P_1^3+ O(\mu), \label{C1 system 5} \\
Q_{2Z} & = P_2  + \tfrac{2}{3}(1+\delta)P_1 + O(\mu), \label{C1 system 6} \\
P_{1Z} & = Q_2+ O(\mu), \label{C1 system 7} \\
P_{2Z} & = Q_1 + \tfrac{2}{3}(1+\delta)Q_2+ O(\mu), \label{C1 system 8}
\end{align}
which is of course reversible with reverser $S:(Q,P) \mapsto (-Q,P)$. Suppose $d_1>0$.
In the limit $(\mu,\check{\kappa}) \rightarrow 0$
equations \eqref{C1 system 5}--\eqref{C1 system 8} are equivalent to the single
fourth-order ordinary differential equation
$$
\partial_Z^4 u-2(1+\delta)\partial_Z^2u+u -u^3 =0
$$
for the variable $u=2d_1^{1/2}P_1$.

Existence theories for homoclinic solutions to the equation
\begin{equation}
\ddddot{u}  - 2(1+\delta) \ddot{u} +u - u^m=0. \label{Generalised BCT equation}
\end{equation}
for a fixed natural number $m \geq 2$  (which is a travelling-wave version of the generalised Kawahara equation)
are given in Theorems \ref{BCT 1} and \ref{BCT 2} below. These theorems are generalisations of
results given by Buffoni, Champneys \& Toland \cite{BuffoniChampneysToland96} (see also Devaney \cite{Devaney76a})
for the special case $m=2$; a full discussion of their generalisation to $m \geq 2$ is given by Ahmad \cite{Ahmad16}.

\begin{theorem} \label{BCT 1}
Suppose that $\delta \geq 0$.
\begin{itemize}
\item[(i)]
Suppose that $m$ is even. Equation \eqref{Generalised BCT equation} has precisely one homoclinic solution $h$ (up to translations). This
solution is positive and symmetric, and monotone increasing to the left, monotone decreasing to the right of its point of symmetry.
\item[(ii)]
Suppose that $m$ is odd. Equation \eqref{Generalised BCT equation} has precisely two homoclinic solutions
$\pm h$, where $h$ is symmetric, and monotone increasing to the left, monotone decreasing to the right of its point of symmetry.
\end{itemize}
In both cases the homoclinic solutions are \underline{transverse}, that is, the stable and unstable manifolds of the zero equilibrium
intersect transversally with respect to the zero level surface of the Hamiltonian at their point of symmetry.
\end{theorem}

\begin{figure}
\centering

\includegraphics[scale=0.7]{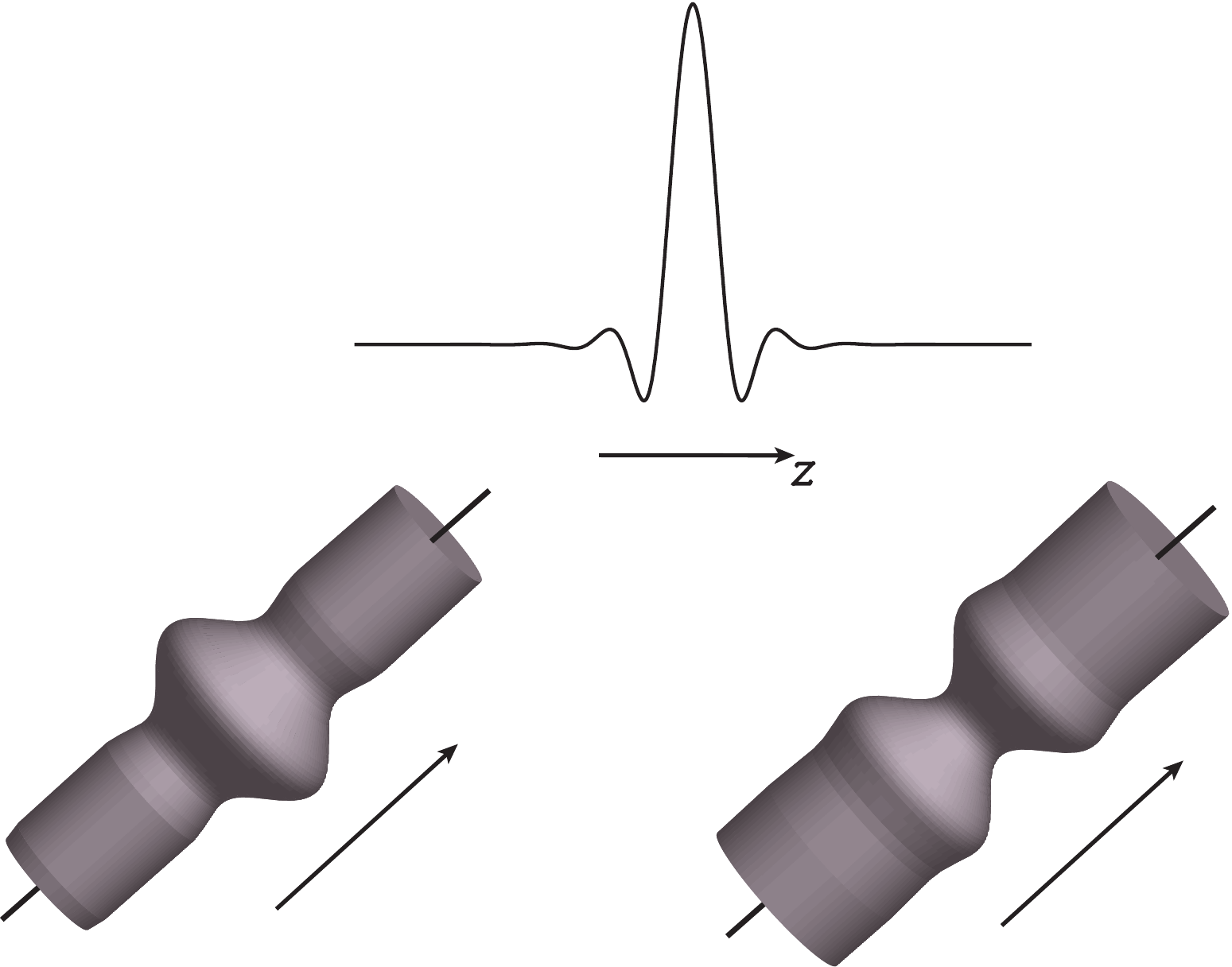}
{\it
\caption{A solitary wave of elevation (left) and depression (right) generated by a `primary'
homoclinic solution (top) in region II} \label{C1 waves - primary}}
\vspace{5mm}
\includegraphics[scale=0.7]{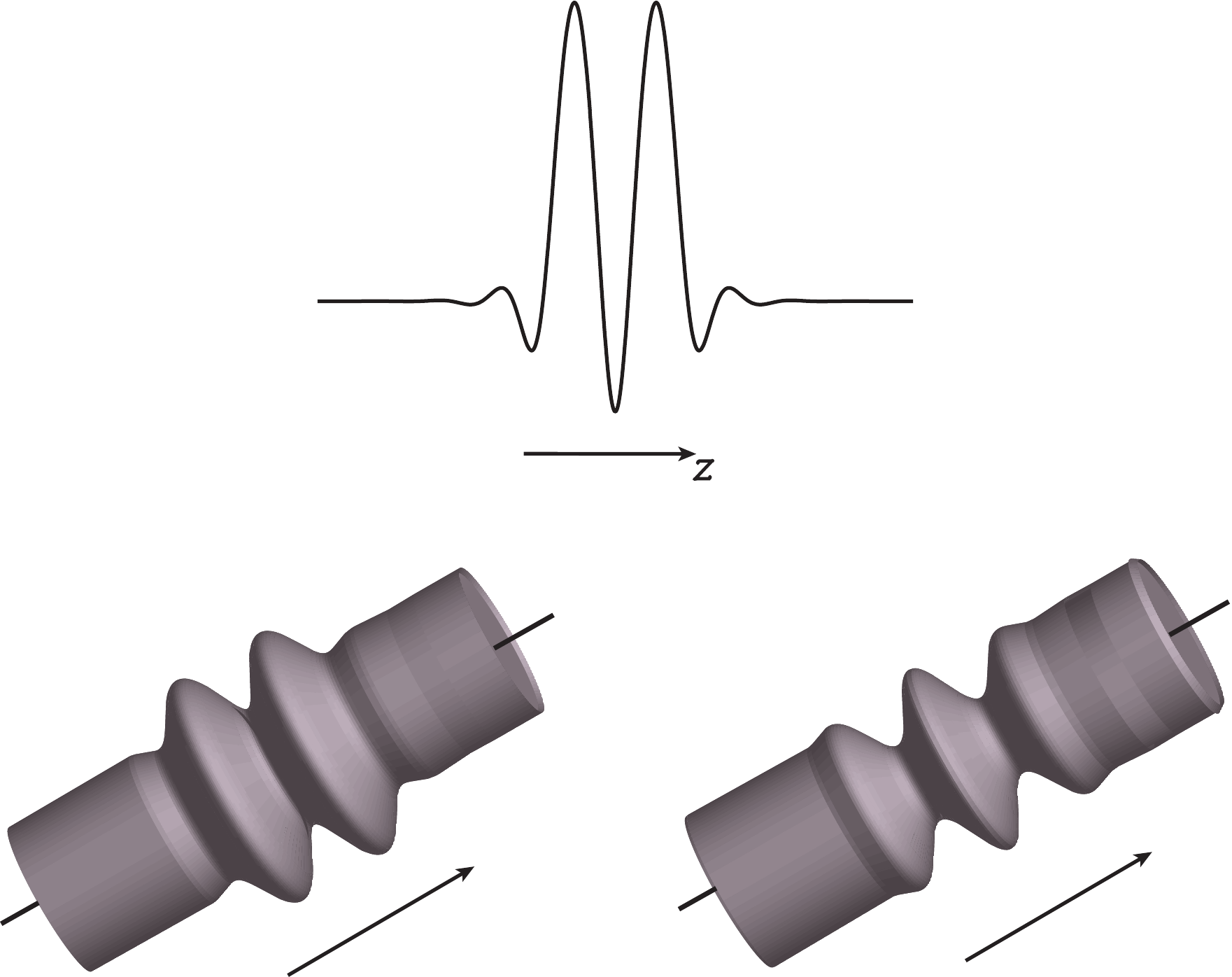}
{\it
\caption{A solitary wave of elevation (left) and depression (right) generated by a `${\bf 2}(2)$'
homoclinic solution (top) in region II}  \label{C1 waves -2(2)}}
\end{figure}

\begin{theorem} \label{BCT 2}
The \underline{primary} homoclinic solutions found in the previous theorem persist
(as small, uniform perturbations of their limits at $\delta=0$) for small, negative values of $\delta$.

Furthermore, each primary homoclinic solution $h$ in the region $\delta<0$ generates a family of transverse \underline{multipulse}
homoclinic solutions which resemble multiple copies of $h$ `glued' together with small oscillations in between.
More precisely, for each all natural numbers $\ell_1, \ldots, \ell_{n-1}$ with $n=1, 2, \ldots$ there exists a homoclinic
solution ${\bf n}(\ell_1,\ldots,\ell_{n-1})$ associated with $h$ which
\begin{itemize}
\item[(i)]
has $n$ local extrema at $t_1$, \ldots, $t_n$,
\item[(ii)]
oscillates $\left\lfloor\frac{\ell_k}{2}\right\rfloor$ times and
has $2\left\lfloor\frac{\ell_k-1}{2}\right\rfloor$ extrema in each 
interval $(t_k,t_{k+1})$,
\item[(iii)]
oscillates infinitely often in the intervals $(-\infty,t_1)$ and $(t_n,\infty)$.
\end{itemize}
\end{theorem}

Theorem \ref{BCT 1}(i) also applies to \eqref{C1 system 1}--\eqref{C1 system 4}
for small, positive values of $\mu$, while Theorem \ref{BCT 1}(ii) applies to \eqref{C1 system 5}--\eqref{C1 system 8}
for small, positive values of $\mu$ and small, values of $\check{\kappa}$ (that is, small values of $\kappa$); the qualitative statements
apply to the variable $\check{c}_1 P_1$ or $\check{d}_1^{1/2}P_1$. The homoclinic orbits at $\mu=0$
(and $\check{\kappa}=0$) are transverse and therefore persist (as small, uniform perturbations
of their limits) for small, positive values of $\mu$ (and small values of $\check{\kappa}$). Similarly, Theorem
\ref{BCT 2} applies to any of these persistent primary homoclinic orbits.

Altogether we have established the existence of a primary and accompanying multipulse family of solitary waves of depression
for $m_1^\prime(1) < \frac{8}{3}$ and elevation for $m_1^\prime(1) > \frac{8}{3}$;
the corresponding ferrofluid surface $\{r=1+\eta(z)\}$ is obtained from the homoclinic
solution of \eqref{Scaled reduced system 1}, \eqref{Scaled reduced system 2}
by the formula
$$\eta(z)=\tfrac{1}{2}\mu^4P_1(\mu z) + O(\mu^5).$$
Furthermore, two multipulse families of solitary waves exist for small values of $m_1^\prime(1)-\frac{8}{3}$ provided that
$m_1^{\prime\prime}(1) \neq \frac{1264}{75}$; one consists of waves of depression, the other of waves of elevation.
The corresponding ferrofluid surface $\{r=1+\eta(z)\}$ is obtained from a homoclinic
solution of \eqref{Scaled reduced system 3}, \eqref{Scaled reduced system 4} by the formula
$$\eta(z)=\tfrac{1}{2}\mu^2P_1(\mu z) + O(\mu^3).$$

\subsection{Homoclinic bifurcation at $C_2$} \label{Region III}

At each point of the curve $C_2$ in Figure \ref{Bifurcation curves} two pairs of purely imaginary eigenvalues become complex by colliding at non-zero points
$\pm \i s$ on the imaginary axis and forming two Jordan chains of length 2. This resonance is associated with the bifurcation of a branch of homoclinic solutions into the region with complex eigenvalues (the parameter regime marked III in Figure \ref{summary}). Let us therefore choose
$$\beta_0=\frac{1}{2}\left(1-\frac{I_0(s)I_2(s)}{I_1(s)^2}\right),
\qquad
\gamma_0 =\frac{1}{2}s^2\left(-1+\frac{I_0(s)^2}{I_1(s)^2}\right)
$$
(so that $\alpha_0=\gamma_0-\beta_0)$ and introduce a bifurcation parameter $\mu$ by writing 
$(\varepsilon_1,\varepsilon_2) = (0,\mu)$, where $0 < \mu \ll 1$.

The six-dimensional centre subspace of $K$ is spanned by the generalised eigenvectors
$$
e_1 = \begin{pmatrix} 0 \\ 0 \\ 1 \\ 0 \end{pmatrix}, \quad
e_2 =\begin{pmatrix} \gamma_0^{-1} \\ 0 \\ 0 \\ 1-2\gamma_0^{-1} \end{pmatrix}, \quad
e, \quad \bar{e}, \quad f, \quad \bar{f},
$$
where
$$
e = \begin{pmatrix} I_1(s) \\ \i s \beta_0 I_1(s)-\i I_2(s) \\ - \i I_0(sr) \\ s I_0(sr)-2I_1(s) \end{pmatrix}, \quad
f = \begin{pmatrix} -\i I_0(s)+\frac{\i}{s} I_1(s) \\
\beta_0 I_0(s) - \frac{2}{s}I_2(s) - I_3(s) \\
-r I_1(s r) \\
-\i I_0(sr)-\i r sI_1(sr)+2i I_0(s) - \frac{2\i}{s} I_1(s)
\end{pmatrix} -\frac{\i\tau_2}{2\tau_1}e
$$
and
\begin{align*}
\tau_1 & = 2 I_0(s)^2 - s\frac{I_0(s)^3}{I_1(s)}+sI_0(s)I_1(s)-I_1(s)^2, \\
\tau_2 & = -\frac{1}{3}\left(\frac{2}{s}(-3+s^2)I_0(s)^2-3s\frac{I_0(s)^4}{I_1(s)^2} + 9 \frac{I_0(s)^3}{I_1(s)}-5I_0(s)I_1(s)+\frac{1}{s}(5+s^2)I_1(s)^2\right);
\end{align*}
note that
$Ke_1=0$, $Ke_2=e_1$, $(K-\i sI) e=0$, $(K-\i sI)f =e$,
$$\Omega(e_1,e_2)=\tfrac{1}{2}-2\gamma_0^{-1}, \quad
\Omega(e,\bar{f})=\tau_1, \quad
\Omega(\bar{e},f)=\tau_1$$
and the symplectic product of any other combination of the vectors $e_1$ $e_2$, $e$, $f$, $\bar{e}$, $\bar{f}$ is zero.
Writing
$$w_1=q_0 f_1 + p_0 f_2 + A E + B F + \bar{A}\bar{E} + \bar{B}\bar{F},$$
where
$$f_1=(\tfrac{1}{2}-2\gamma_0^{-1/2})^{-1}e_1, \qquad f_2=(\tfrac{1}{2}-2\gamma_0^{-1})^{-1/2}e_2, \qquad E=\tau_1^{-1/2}e, \qquad F = \tau_1^{-1/2}f,$$
we therefore find that $q_0$, $p_0$, $A$ and $B$ are canonical coordinates for the reduced Hamiltonian system,
which has the cyclic variable $q_0$ and reverser $S:(q_0,p_0,A,B) \mapsto (-q_0,p_0,\bar{A},-\bar{B})$; with a slight abuse
of notation we abbreviate $\tilde{H}^{\varepsilon}|_{(\varepsilon_1,\varepsilon_2)=(0,\mu)}$ to $\tilde{H}^\mu$.

The usual normal-form theory for the two-dimensional system with Hamiltonian $\tilde{H}^\mu(A,B,\bar{A},\bar{B},0)$ asserts that,
after a canonical change of variables,
$$\tilde{H}^\mu(A,B,\bar{A},\bar{B},0) \!=\! \i s (A\bar{B}-\bar{A}B)+|B|^2
 + \tilde{H}_\mathrm{NF}^0(|A|^2,\i(A\bar{B}-\bar{A}B),\mu)
+ O(|(A,B)|^2|(\mu,A,B)|^{n_0}),$$
where $\tilde{H}_\mathrm{NF}^0$ is a real polynomial function of its arguments which satisfies
\[\tilde{H}_\mathrm{NF}^0(|A|^2,\i(A\bar{B}-\bar{A}B),\mu)=O(|(A,B)|^2|(\mu,A,B)|).\]
It follows that, after a canonical change of variables,
$$
\tilde{H}^\mu(A,B,\bar{A},\bar{B},p_0) = \i s (A\bar{B}-\bar{A}B)+|B|^2 + \tfrac{1}{2}p_0^2 + \tilde{H}^\mu_\mathrm{nl}(A,B,\bar{A},\bar{B},p_0)
$$
with
\begin{align*}
\tilde{H}^\mu_\mathrm{nl}(A,B,\bar{A},\bar{B},p_0)
& = \tilde{H}_\mathrm{NF}(|A|^2,\i(A\bar{B}-\bar{A}B),p_0,\mu)\\
& \qquad\qquad\mbox{} + \tilde{H}_\mathrm{r}(A,B,\bar{A},\bar{B},p_0,\mu)
+ O(|(A,B,p_0)|^2|(\mu,A,B,p_0)|^{n_0});
\end{align*}
here $\tilde{H}_\mathrm{NF}$ is a real polynomial function of its arguments which satisfies
\[\tilde{H}_\mathrm{NF}(|A|^2,\i(A\bar{B}-\bar{A}B),p_0,\mu)=O(|(A,B)|^2|(\mu,A,B,p_0)|)\]
and $\tilde{H}_\mathrm{NF}(|A|^2,\i(A\bar{B}-\bar{A}B),0,\mu)=\tilde{H}_\mathrm{NF}^0(|A|^2,\i(A\bar{B}-\bar{A}B),\mu)$,
and $\tilde{H}_\mathrm{r}(A,B,\bar{A},\bar{B},p_0,\mu)$ is an affine function of its first four arguments which satisfies
\[\tilde{H}_\mathrm{r}(|A|^2,\i(A\bar{B}-\bar{A}B),p_0,\mu)=O(|(A,B,p_0)||p_0||(\mu,A,B,p_0)|)\]
Note that
\begin{align*}
P^\mu(A,B,\bar{A},\bar{B},p_0) & = \partial_{\bar{B}} \tilde{H}_\mathrm{nl}^\mu(A,B,\bar{A},\bar{B},p_0)E + \partial_B \tilde{H}_\mathrm{nl}^\mu(A,B,\bar{A},\bar{B},p_0)\bar{E}\\
& \qquad\qquad\mbox{}
-\partial_{\bar{A}} \tilde{H}_\mathrm{nl}^\mu(A,B,\bar{A},\bar{B},p_0)F - \partial_A \tilde{H}_\mathrm{nl}^\mu(A,B,\bar{A},\bar{B},p_0)\bar{F} 
\\
& \qquad\qquad \mbox{}
+\partial_{p_0} \tilde{H}_\mathrm{nl}^\mu(A,B,\bar{A},\bar{B},p_0)f_1.
\end{align*}

Writing
\begin{align*}
\tilde{H}_2^1(A,B,p_0) & = c_1^1p_0^2+c_2^1|A|^2+c_3^1\i(A\bar{B}-\bar{A}B)+c_4^1p_0A+\bar{c}_4^1p_0\bar{A}
+c_5^1p_0B+\bar{c}_5^1p_0\bar{B},\\
\tilde{H}_3^0(A,B,p_0)&=c_1p_0^3+c_2p_0\vert A\vert^2+c_3\i p_0(A\bar{B}-\bar{A}B)+c_4p_0^2A+\bar{c}_4p_0^2\bar{A}+c_5p_0^2B+\bar{c}_5p_0^2\bar{B}\\
\tilde{H}_4^0(A,B,p_0) & = d_1 p_0^4 +d_2 p_0^2|A|^2 + d_3\i p_0^2(A\bar{B}-\bar{A}B)+d_4 |A|^4
+d_5\i (A\bar{B}-\bar{A}B)|A|^2 \\
& \qquad\mbox{} -d_6(A\bar{B}-\bar{A}B)^2 + d_7 p_0^3 A + \bar{d}_7p_0^3 \bar{A}+ d_8 p_0^3 B +\bar{d}_8p_0^3 \bar{B}, 
\end{align*}
where $\mu^j\tilde{H}^j_k(A,B,p_0)$ denotes
the part of the Taylor expansion of $\tilde{H}^\mu(A,B,p_0)$
which is homogeneous of order $j$ in $\mu$ and $k$ in $(A,B,p_0)$,
one finds that
\small
\begin{align}
d_4 & = \frac{I_1(s)^2}{2\tau_1^2}\!\Bigg(\! \frac{(-2s^2+s^2\beta_0-2sT+4s^2ST-\alpha_0m_1^\prime(1))(-2s^2-2s^2\beta_0-sS+4s^2ST-\alpha_0m_1^\prime(1))}{2(\gamma_0+4s^2\beta_0-2sT)} \nonumber \\
& \hspace{0.8in} \mbox{}-\frac{(s^2\beta_0-4sS+2+\alpha_0m_1^\prime(1))(3sS-\alpha_0m_1^\prime(1))}{\gamma_0-2} \nonumber \\
& \hspace{0.8in} \mbox{}+7s^2-\tfrac{21}{2}s^2\beta_0+\tfrac{3}{2}s^4\beta_0+6sS-6s^3S+4s^3S^2T-2s^2ST-3\alpha_0m_1^\prime(1)-\tfrac{1}{2}\alpha_0m_1^{\prime\prime}(1)\!\!\Bigg), \nonumber \\
c_2^1 &= - \frac{I_1(s)^2}{\tau_1}, \label{Formula for c4}
\end{align}
\normalsize
where
$$S=\frac{I_0(s)}{I_1(s)}, \qquad T=\frac{I_0(2s)}{I_1(2s)}$$
(see Appendix (v)).

The lower-order Hamiltonian system
\begin{align}
A_Z &= \partial_{\bar{B}} \tilde{H}^\mu(A,B,\bar{A},\bar{B},0), \label{IP BG 1}, \\
B_Z &= -\partial_{\bar{A}} \tilde{H}^\mu(A,B,\bar{A},\bar{B},0) \label{IP BG 2}
\end{align}
has been examined in detail by Iooss \& P\'{e}rou\`{e}me \cite{IoossPeroueme93}. The
`truncated normal form' obtained by ignoring the remainder terms in $ \tilde{H}^\mu(A,B,\bar{A},\bar{B},0)$
is conveniently handled using the substitution
$A(z)=\mathrm{e}^{\i s z}a(z)$, $B(z)=\mathrm{e}^{\i s z}b(z)$, which converts it into the system
\begin{align}
\dot{a} & = b + \partial_b\tilde{H}_\mathrm{NF}^0(|a|^2,\i(a\bar{b}-\bar{a}b),\mu), \label{IP reduction 1} \\
\dot{b} & = -\partial_{\bar{a}}\tilde{H}_\mathrm{NF}^0(|a|^2,\i(a\bar{b}-\bar{a}b),\mu). \label{IP reduction 2}
\end{align}
Supposing that the coefficients $c_2^1$ and $d_4$ are respectively negative and positive,
one finds that \eqref{IP reduction 1}, \eqref{IP reduction 2}
admits a real, reversible homoclinic solution $(a_\mathrm{h},b_\mathrm{h})$, which evidently generates a circle
$\{\mathrm{e}^{\i \theta}(a_\mathrm{h},b_\mathrm{h}): \theta \in [0,2\pi)\}$ of further homoclinic solutions, two of
which (those with $\theta=0$ and $\theta=\pi$) are reversible. The corresponding pair of homoclinic solutions
to the original `truncated normal form' are reversible and persist when the remainder terms are
reinstated. A theory of multipulse homoclinic solutions to \eqref{IP BG 1}, \eqref{IP BG 2} has also been
given by Buffoni \& Groves \cite{BuffoniGroves99} (under the same hypotheses on the normal-form coefficients).

\begin{figure}[h]
\centering
\includegraphics[scale=0.5]{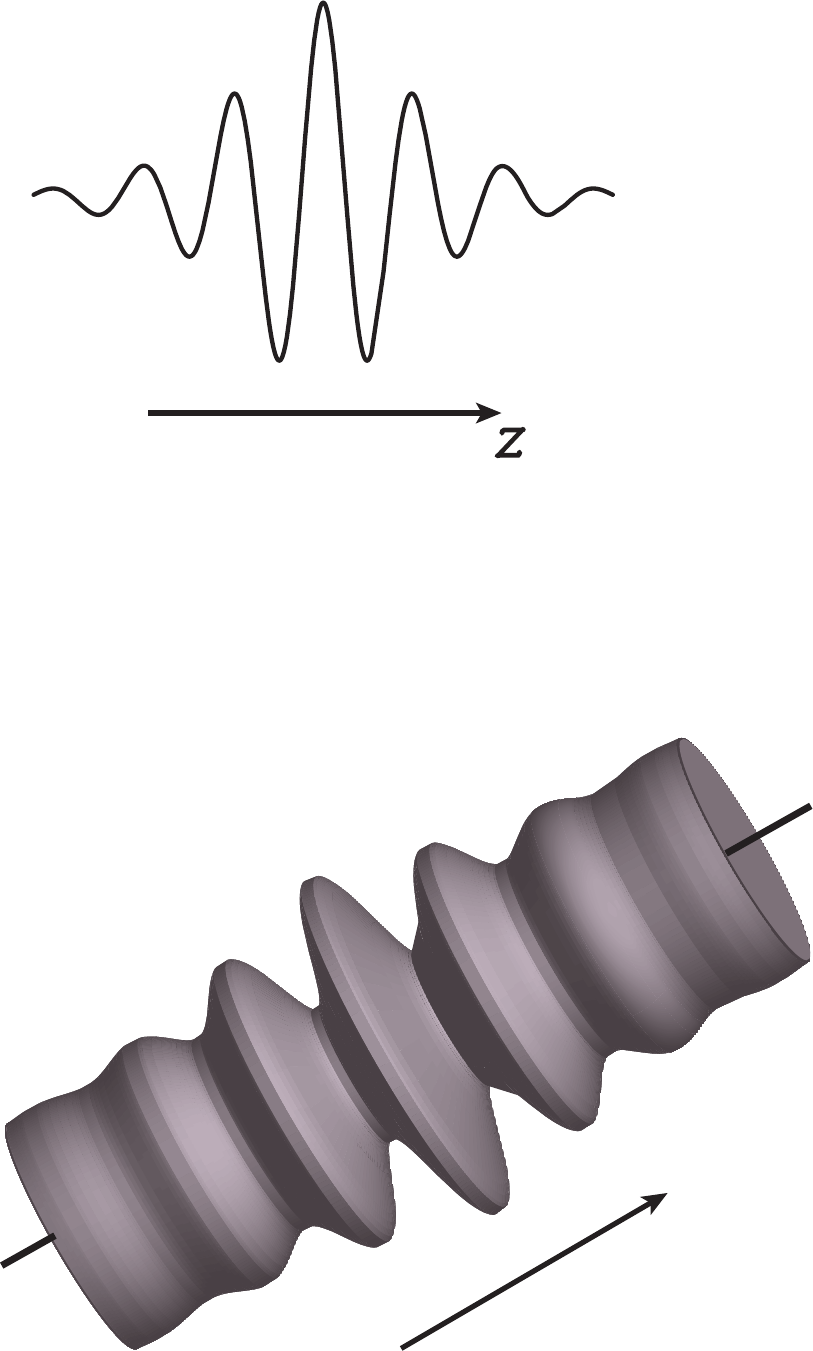}
\hspace{1.5cm}\includegraphics[scale=0.18]{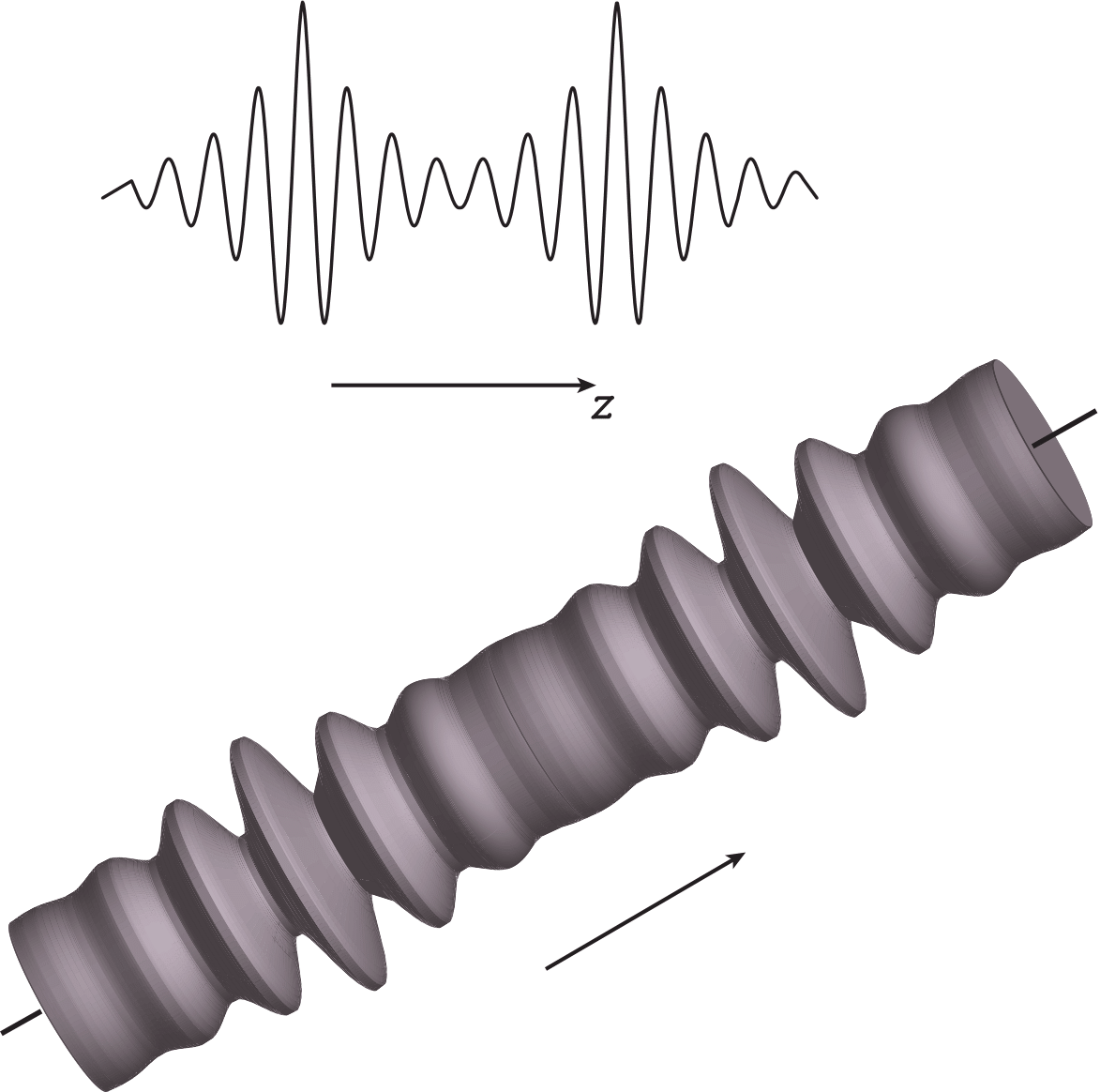}

{\it \caption{Sketches of $\eta=\eta(z)$ and the corresponding symmetric unipulse (left) and
multipulse solitary waves (right) generated by homoclinic solutions in region III}
\label{ESW sketches}}
\end{figure}

\begin{theorem} $ $
\begin{list}{(\roman{count})}{\usecounter{count}}
\item (Iooss \& P\'{e}rou\`{e}me) For each sufficiently small,
positive value of $\mu$ the 
two-degree-of-freedom Hamiltonian system \eqref{IP BG 1}, \eqref{IP BG 2}
has two distinct symmetric homoclinic 
solutions.
\item (Buffoni \& Groves) For each sufficiently small, positive value of $\mu$ the 
two-degree-of-freedom Hamiltonian system \eqref{IP BG 1}, \eqref{IP BG 2}
has an infinite number of geometrically 
distinct homoclinic solutions which generically resemble multiple
copies of one of the homoclinic solutions in part (i).
\end{list}
The homoclinic solutions identified above correspond to envelope solitary waves whose amplitude
is $O(\mu^{1/2})$ and which decay exponentially as $z \rightarrow \pm\infty$;
they are sketched in Figure \ref{ESW sketches}. 
\end{theorem}

\section*{Appendix: Calculation of the normal-form coefficients}
\renewcommand{\theequation}{A.\arabic{equation}}
\setcounter{equation}{0}

The coefficients in the reduced Hamiltonian $\tilde{H}^\varepsilon(w_1)$ are determined using
the equations
\begin{align}
K \tilde{r}(w_1;\varepsilon) - \mathrm{d}_1\tilde{r}[w_1;\varepsilon](Kw_1) & = P^\varepsilon(w_1)
+ \mathrm{d}_1\tilde{r}[w_1;\varepsilon](P^\varepsilon(w_1)) - g_\mathrm{nl}^\varepsilon(w_1+\tilde{r}(w_1;\varepsilon)), \label{Cohomo 1}\\
B_\mathrm{l}\tilde{r}(w;\varepsilon) & = - B_\mathrm{nl}^\varepsilon(w_1+\tilde{r}(w_1;\varepsilon)) \label{Cohomo 2}
\end{align}
to compute the Taylor series of $\tilde{H}^\varepsilon(w_1)$ and $\tilde{r}(w_1;\varepsilon)$ systematically in powers of $(q,p,q_0)$ or $(q,p,p_0)$. Here
$K=\mathrm{d}g^0[0]$ and $g_\mathrm{nl}^\varepsilon = g^\varepsilon - K$ are the linear and nonlinear
parts of $g^\varepsilon$ (with this slight abuse of notation $K$ is given by the explicit formula \eqref{Formula for K}),
and $B_\mathrm{l}$, $B_\mathrm{nl}^\varepsilon$ are the linear and nonlinear parts of the boundary-value operator
$B^\varepsilon: N \rightarrow {\mathbb R}$ defined by the left-hand side of \eqref{BC 2}.
Throughout these calculations we also make use of the identity
$$
\Omega(Ku+g_\mathrm{nl}^\varepsilon(u),v)+(B_\mathrm{l}(u)+B_\mathrm{nl}^\varepsilon(u))\phi^v\vert_{r=1}=\mathrm{d} H^\varepsilon[u](v),
$$
in which $v=(\eta^v,\omega^v,\phi^v,\zeta^v)$.
We denote the parts of $H^\varepsilon(w)$, $B_\mathrm{nl}^\varepsilon(w)$, $g_\mathrm{nl}^\varepsilon(w)$
which are homogeneous of order $m$ in $\varepsilon$ and $n$ in $w$ by
$\varepsilon^m H_n^m(w)$, $\varepsilon^m B_{\mathrm{nl},n}^m(w)$, $\varepsilon^m g_{\mathrm{nl},n}^m(w)$,
and the part of $\tilde{r}(w_1;\varepsilon)$ which is homogeneous of order $m$ in $\varepsilon$ and $n$ in $w_1$
by $\tilde{r}_n^m(w_1;\varepsilon)$; the notation is modified in the natural fashion when $\varepsilon$ is replaced
by a more specific parameterisation. Finally, arbitrary constants arising from solving differential equations are denoted
by $a_i$.

\subsubsection*{Homoclinic bifurcation at $C_4$}\indent

(i) Write
$$
\tilde{r}_m^n(w_1;\mu)=\sum_{h+i+j=m}\mu^n\tilde{r}_{hij0}^nq^hp^iq_0^j
$$
and consider the $q^2$ and $\mu q$ components of \eqref{Cohomo 1}, \eqref{Cohomo 2}, namely
\begin{align}\label{q2c4}
q^2:\ 
&\begin{cases}
K\tilde{r}_{2000}^0=-3c_1f_3-c_2f_1-g_{\mathrm{nl},2}^0(f_2,f_2), \\
B_\mathrm{l}\tilde{r}_{2000}^0=-B_{\mathrm{nl},2}^0(f_2,f_2),
\end{cases} \\[1mm]
\mu q:\ 
&\begin{cases}
K\tilde{r}_{1000}^1=-2c_1^1f_3-c_2^1f_1-g_{\mathrm{nl},1}^1(f_2),\\
B_\mathrm{l}\tilde{r}_{1000}^1=0.
\end{cases} \nonumber
\end{align}
Using these equations we find that
\begin{align*}
c_1&=H_3^0(f_2,f_2,f_2)+2H_2^0(\tilde{r}_{2000}^0,f_2)\\
&=H_3^0(f_2,f_2,f_2)+\Omega(K\tilde{r}_{2000}^0,f_2)+B_\mathrm{l}\tilde{r}_{2000}^0\phi^{f_2}\vert_{r=1}\\
&=H_3^0(f_2,f_2,f_2)+3c_1-\Omega(g_{\mathrm{nl},2}^0(f_2,f_2),f_2)\\
&=-2H_3^0(f_2,f_2,f_2)+3c_1,
\end{align*}
which implies that
$$
c_1=H_3^0(f_2,f_2,f_2)=\tfrac{1}{6}\left(\beta_0-\tfrac{1}{4}\right)^{-3/2}(\alpha_0 m_1^\prime(1)-6),
$$
and
$$
c_1^1=H_2^1(f_2,f_2)+2H_2^0(\tilde{r}_{1000}^1,f_2)
=H_2^1(f_2,f_2)+2c_1^1-\Omega(g_{\mathrm{nl},1}^1(f_2),f_2)
=-H_2^1(f_2,f_2)+2c_1^1,
$$
which implies that
$$c_1^1=H_2^1(f_2,f_2)=-\tfrac{1}{2}\left(\beta_0-\tfrac{1}{4}\right)^{-1}.$$

(ii) Write
$$
\tilde{r}_m^{0,0}(w_1;\mu,\kappa)=\sum_{h+i+j=m}\tilde{r}_{hij0}^{0,0} q^hp^iq_0^j
$$
and consider the $q^3$ component of \eqref{Cohomo 1}, \eqref{Cohomo 2}, namely
\begin{equation}\label{q3c4}
q^3: \ 
\begin{cases}
K\tilde{r}_{3000}^{0,0}=-4e_1f_3-e_2f_1-g_{\mathrm{nl},3}^{0,0}(f_2,f_2,f_2)-2g_{\mathrm{nl},2}^{0,0}(f_2,\tilde{r}_{2000}^{0,0}), \\
B_\mathrm{l}\tilde{r}_{3000}^{0,0}=-B_{\mathrm{nl},3}^{0,0}(f_2,f_2,f_2)-2B_{\mathrm{nl},2}^{0,0}(f_2,\tilde{r}_{2000}^{0,0}).
\end{cases}
\end{equation}
The coefficient $d_1$ can be expressed as
\begin{equation}\label{d1c4}
d_1=H_4^{0,0}(f_2,f_2,f_2,f_2)+3H_3^{0,0}(f_2,f_2,\tilde{r}_{2000}^{0,0})+2H_2^{0,0}(f_2,\tilde{r}_{3000}^{0,0})+H_2^{0,0}(\tilde{r}_{2000}^{0,0}, \tilde{r}_{2000}^{0,0}),
\end{equation}
where
\begin{align*}
2H_2^{0,0}(f_2,\tilde{r}_{3000}^{0,0})&=\Omega(K\tilde{r}_{3000}^{0,0},f_2)+B_\mathrm{l}\tilde{r}_{3000}^{0,0}\phi^{f_2}\vert_{r=1}\\
&=4d_1-\Omega(g_{\mathrm{nl},3}^{0,0}(f_2,f_2,f_2),f_2)-2\Omega(g_{\mathrm{nl},2}^{0,0}(\tilde{r}_{2000}^{0,0},f_2),f_2)+B_\mathrm{l}\tilde{r}_{3000}^{0,0}\phi^{f_2}\vert_{r=1}\\
&=4d_1-4H_4^{0,0}(f_2,f_2,f_2,f_2)-6H_3^{0,0}(f_2,f_2,\tilde{r}_{3000}^{0,0})\\
&\qquad \mbox{}+\big(B_\mathrm{l}\tilde{r}_{3000}^{0,0}+2B_{\mathrm{nl},2}^{0,0}(f_2,\tilde{r}_{2000}^{0,0})
+B_{\mathrm{nl},3}(f_2,f_2,f_2)\big)\phi^{f_2}\vert_{r=1},
\end{align*}
and we find from the boundary condition in \eqref{q3c4} that the sum inside the parentheses vanishes.
From \eqref{q2c4} we find that
$$
\tilde{r}_{2000}^{0,0}=\big( 6(\beta_0-\tfrac{1}{4})^{-1/2}-c_2\big)f_2+a_1f_1,
$$
and it follows from
\eqref{d1c4} that
$$
d_1=H_4^{0,0}(f_2,f_2,f_2,f_2)+H_3^{0,0}(f_2,f_2,\tilde{r}_{2000}^{0,0})-\tfrac{1}{3}H_2^{0,0}(\tilde{r}_{2000}^{0,0},\tilde{r}_{2000}^{0,0})
=\tfrac{1}{24}(\beta_0-\tfrac{1}{4})^{-2}(12-\alpha_0 m_1^{\prime\prime}(1)).
$$

\subsubsection*{Homoclinic bifurcation at $C_1$}\indent

(iii) Write
$$
\tilde{r}_m^{n_1,n_2}(w_1;\mu_1,\mu_2)=\sum_{h+i+j+k+l=m}\mu_1^{n_1}\mu_2^{n_2}
\tilde{r}_{hijk0l}^{m_1,m_2}q_1^hp_1^iq_2^jp_2^kp_0^l
$$
and consider the $p_1^2$, $\mu_2 p_1$ and $\mu_1 p_1$ components of \eqref{Cohomo 1}, \eqref{Cohomo 2}, namely
\begin{align}\label{p12c1}
p_1^2:\ &\begin{cases}
K\tilde{r}_{020000}^{0,0}=3c_1f_5-2c_5f_3+c_2f_1-g_{\mathrm{nl},2}^{0,0}(f_2,f_2),\\
B_\mathrm{l}\tilde{r}_{020000}^{0,0}=-B_{\mathrm{nl},2}^{0,0}(f_2,f_2),
\end{cases} \\[1mm]
\mu_2p_1:\ &\begin{cases}
K\tilde{r}_{010000}^{0,1}=2c_1^{0,1}f_5-2c_4^{0,1}f_3+c_2^{0,1}f_1-g_{\mathrm{nl},1}^{0,1}(f_2),\\
B_\mathrm{l}\tilde{r}_{010000}^{0,1}=0,
\end{cases} \nonumber \\[1mm]
\label{epsilonp1c1}
\mu_1p_1:\ & \begin{cases}
K\tilde{r}_{010000}^{1,0}=2c_1^{1,0}f_5-2c_4^{0,1}f_3+c_2^{1,0}f_1-g_{\mathrm{nl},1}^{1,0}(f_2), \\
B_\mathrm{l}\tilde{r}_{010000}^{1,0}=-B_{\mathrm{nl},1}^{1,0}(f_2).
\end{cases}
\end{align}
Using the method described in part (i) above, we find from these equations that
\begin{align*}
c_1&=H_3^{0,0}(f_2,f_2,f_2)=48\sqrt{6}(3 m_1^\prime(1)-8), \\
c_1^{0,1}&=H_2^{0,1}(f_2,f_2)=-48, \\
c_1^{1,0}&=H_2^{1,0}(f_2,f_2)=0.
\end{align*}

Combining
\begin{equation}\label{epsilon2q2}
\mu_1q_2:\ \begin{cases}
K\tilde{r}_{001000}^{1,0}-\tilde{r}_{010000}^{1,0}=-2c_4^{1,0}f_4-g_{\mathrm{nl},1}^{1,0}(f_3),\\
B_\mathrm{l}\tilde{r}_{001000}^{1,0}=-B_{\mathrm{nl},1}^{1,0}(f_3),
\end{cases}
\end{equation}
with
$$
\tilde{r}_{010000}^{1,0}=-2c_4^{1,0}f_4+c_2^{1,0}f_2+a_2f_1,
$$
which is obtained from \eqref{epsilonp1c1}, one finds by the usual argument that
$$
c_4^{1,0}=\frac{1}{3}H_2^{1,0}(f_3,f_3)=-16.
$$
Similarly, combining
\begin{equation}\label{epsilon12p1}
\mu_1^2p_1:\ \begin{cases}
K\tilde{r}_{010000}^{2,0}=2c_1^{2,0}f_5-2c_4^{2,0}f_3+c_2^{2,0}f_1-2c_4^{1,0}\tilde{r}_{001000}^{1,0}, \\
B_\mathrm{l}\tilde{r}_{010000}^{2,0}=0
\end{cases}
\end{equation}
with
$$
\tilde{r}_{001000}^{1,0}=-4c_4^{1,0}f_5+c_2^{1,0}f_3+a_2f_2+a_3f_1+
\left(\begin{array}{c}
0\\
4\sqrt{6}\\
0\\
0
\end{array}
\right),
$$
which is obtained from \eqref{epsilon2q2}, yields
\begin{align*}
c_1^{2,0}&=H_2^{2,0}(f_2,f_2)+2H_2^{0,0}(\tilde{r}_{010000}^{2,0},f_2)+2H_2^{1,0}(\tilde{r}_{010000}^{1,0},f_2)+H_2^{0,0}(\tilde{r}_{010000}^{1,0},\tilde{r}_{010000}^{1,0})\\
&=2c_1^{2,0}-512,
\end{align*}
so that
$ c_1^{2,0}=512$. 

(iv) Write
$$
\tilde{r}_m^{0,0,0}(w_1;\mu_1,\mu_2,\kappa)=\sum_{h+i+j+k+l=m}\tilde{r}_{hijik0l}^{0,0,0}q_1^hp_1^iq_2^jp_2^kp_0^l
$$
and note that
$$
d_1=H_4^{0,0,0}(f_2,f_2,f_2,f_2)+3H_3^{0,0,0}(f_2,f_2,\tilde{r}_{020000}^{0,0,0})+2H_2^{0,0,0}(f_2,\tilde{r}_{030000}^{0,0,0})
+H_2^{0,0,0}(\tilde{r}_{020000}^{0,0,0},\tilde{r}_{020000}^{0,0,0}).
$$
Since
$$
2H_2^{0,0,0}(f_2,\tilde{r}_{030000}^{0,0,0})
=4d_1-4H_4^{0,0,0}(f_2,f_2,f_2,f_2)-6H_3^{0,0,0}(f_2,f_2,\tilde{r}_{020000}^{0,0,0})-2c_5\Omega(\tilde{r}_{011000}^{0,0,0},f_2),
$$
where we have used
$$
p_1^3:\ \begin{cases}
K\tilde{r}_{030000}^{0,0,0}=4d_1f_5+d_2f_1-2d_7f_3-g_{\mathrm{nl},3}^{0,0,0}(f_2,f_2,f_2) -2g_{\mathrm{nl},2}^{0,0,0}(f_2,\tilde{r}_{020000}^{0,0,0})-2c_5\tilde{r}_{011000}^{0,0,0},\\
B_\mathrm{l}\tilde{r}_{030000}^{0,0,0}=-B_{\mathrm{nl},3}^{0,0,0}(f_2,f_2,f_2)-2B_{\mathrm{nl},2}^{0,0,0}(f_2,\tilde{r}_{020000}^{0,0,0}),\end{cases}
$$
it follows that
\begin{equation}\label{d4formula}
3d_1=3H_4^{0,0,0}(f_2,f_2,f_2,f_2)+3H_3^{0,0,0}(f_2,f_2,\tilde{r}_{020000}^{0,0,0})-H_2^{0,0,0}(\tilde{r}_{020000}^{0,0,0},\tilde{r}_{020000}^{0,0,0})+2c_5\Omega(\tilde{r}_{011000}^{0,0,0},f_2).
\end{equation}
In order to compute $d_1$ it is therefore necessary to compute
$\tilde{r}_{020000}^{0,0,0},\ \tilde{r}_{011000}^{0,0,0}$ and $c_5$.

From \eqref{p12c1} one finds that
$$
\tilde{r}_{020000}^{0,0,0}=-2c_5f_4+(c_2+6\sqrt{6})f_2+a_4f_1,
$$
and
$$
p_1q_2:\ \begin{cases}
K\tilde{r}_{011000}^{0,0,0}-2\tilde{r}_{020000}^{0,0,0}=-2c_5f_4-2g_{\mathrm{nl},2}^{0,0,0}(f_2,f_3), \\
B_\mathrm{l}\tilde{r}_{011000}^{0,0,0}=-2B_{\mathrm{nl},2}^{0,0,0}(f_2,f_3)
\end{cases}
$$
yields
$$
\tilde{r}_{011000}^{0,0,0}=-6c_5f_5+(2c_2+12\sqrt{6})f_3+a_4f_2+a_5f_1+
\left(\begin{array}{c}
0\\
72\\
48r^2\\
0
\end{array}\right).
$$
Furthermore,
$$
q_2^2:\ \begin{cases}
K\tilde{r}_{002000}^{0,0,0}-\tilde{r}_{011000}^{0,0,0}=c_5f_5+c_6f_1-g_{\mathrm{nl},2}^{0,0,0}(f_3,f_3),\\
 B_\mathrm{l}\tilde{r}_{002000}^{0,0,0}=-B_{\mathrm{nl},2}^{0,0,0}(f_3,f_3),
\end{cases}
$$
and
$$
p_1p_2:\ \begin{cases}
K\tilde{r}_{010100}^{0,0,0}-\tilde{r}_{011000}^{0,0,0}=-4c_5f_5-2c_6-2g_{\mathrm{nl},2}^{0,0,0}(f_2,f_4),\\
B_\mathrm{l}\tilde{r}_{010100}^{0,0,0}=-2B_{\mathrm{nl},2}^{0,0,0}(f_2,f_4)
\end{cases}
$$
yield
\begin{align*}
\tilde{r}_{002000}^{000}&=-5c_5f_6+(2c_2+12\sqrt{6})f_4+c_6f_2+a_4f_3+a_5f_2+a_6f_1+
\left(\begin{array}{c}
54\\
0\\
0\\
-108+144r^2
\end{array}\right), \\
\tilde{r}_{010100}^{0,0,0}&=-10c_5f_6+(2c_2+12\sqrt{6})f_4-2c_6f_2+a_4f_3+a_5f_2+a_7f_1+
\left(\begin{array}{c}
39\\
0\\
0\\
-48(1+r^2)
\end{array}\right),
\end{align*}
and using these results we find that the solvability condition for
$$
q_2p_2:\ \begin{cases}
K\tilde{r}_{001100}^{0,0,0}-2\tilde{r}_{002000}^{0,0,0}-\tilde{r}_{010100}^{0,0,0}=-2g_{\mathrm{nl},2}^{0,0,0}(f_3,f_4),\\
B_\mathrm{l}\tilde{r}_{001100}^{0,0,0}=-2B_{\mathrm{nl},2}^{0,0,0}(f_3,f_4)
\end{cases}
$$
is $c_5=-\frac{144\sqrt{6}}{5}$.
Inserting these expressions for $\tilde{r}_{020000}^{0,0,0},\ \tilde{r}_{011000}^{0,0,0}$ and $c_5$ into \eqref{d4formula}, we obtain
$$
d_1=864\left(\frac{1264}{75}-m_1^\prime(1))\right).
$$

\subsubsection*{Homoclinic bifurcation at $C_2$}\indent

(v) Here we write
$$
\tilde{r}_m^n(w_1;\mu)=\hspace{-4mm}\sum_{h+i+j+k+l=m}\tilde{r}_{hijk0l}^n\mu^n A^hB^i\bar{A}^j\bar{B}^kp_0^l.
$$
The coefficient $c_2^1$ is found from
\begin{equation}\label{mu_a}
\mu A:\ \begin{cases}
(K-\i s I)\tilde{r}_{100000}^1=c_3^1\i E-c_2^1F+c_4^1f_1-g_{\mathrm{nl},1}^1(E),\\
B_\mathrm{l}\tilde{r}_{100000}^1=0.
\end{cases}
\end{equation}
Noting that
$$
\Omega(K\tilde{r}_{100000}^1,\bar{E})=2H_2^0(\tilde{r}_{100000}^0,\bar{E})=\Omega(K\bar{E},\tilde{r}_{100000}^0),
$$
we find from \eqref{mu_a} that
$$
c_2^1=-\Omega(\tilde{r}_{100000}^1,(K+\i s I)\bar{E})+\Omega(g_{\mathrm{nl},1}^1(E),\bar{E})=2H_2^1(E,\bar{E})=-\frac{I_1(s)^2}{\tau_1}.
$$ 

Finally, to compute $d_4$ we consider
$$
A|A|^2:\ \begin{cases}
(K-\i s I)\tilde{r}_{201000}^0=\i d_5E-2d_4F-3g_{\mathrm{nl},3}^0(E,E,\bar{E})-2g_{\mathrm{nl},2}^0(\bar{E},\tilde{r}_{200000}^0)-2g_{\mathrm{nl},2}^0(E,\tilde{r}_{101000}^0),\\
B_\mathrm{l}\tilde{r}_{201000}^0=-3B_{\mathrm{nl},3}^0(E,E,\bar{E})-2B_{\mathrm{nl},2}^0(\bar{E},\tilde{r}_{200000}^0)-2B_{\mathrm{nl},2}^0(E,\tilde{r}_{101000}^0).
\end{cases}
$$
Taking the symplectic product with $\bar{E}$ and simplifying in the usual fashion, we find that
\begin{equation}\label{c41}
d_4=6H_4(E,E,\bar{E},\bar{E})+3H_3^0(\tilde{r}_{200000}^0,\bar{E},\bar{E})+3H_3^0(\tilde{r}_{101000}^0,E,\bar{E}),
\end{equation}
where $\tilde{r}_{200000}^0$ and $\tilde{r}_{101000}^0$ are obtained from
$$
A^2:\ \begin{cases}
(K-2\i sI)\tilde{r}_{200000}^0=-g_{\mathrm{nl},2}^0(E,E),\\
B_\mathrm{l}\tilde{r}_{200000}^0=-B_{\mathrm{nl},2}^0(E,E)
\end{cases}
$$
and
$$
|A|^2:\ \begin{cases}
K(\tilde{r}_{101000}^0-c_2f_2)=-2g_{\mathrm{nl},2}^0(E,\bar{E}),\\
B_\mathrm{l}(\tilde{r}_{101000}^0-c_2f_2)=0,
\end{cases}
$$
where
$$
c_2=6H_3^0(E,\bar{E},f_2)
$$
because of
$$
p_0A:\ \begin{cases}
(K-\i s I)\tilde{r}_{100001}^0=\i c_3E-c_2F+2c_4f_1-2g_{\mathrm{nl},2}^0(E,f_2),\\
B_\mathrm{l}\tilde{r}_{100001}^0=-2B_{\mathrm{nl},2}^0(E,f_2)
\end{cases}
$$
(note that $\tilde{r}_{101000}^0$ is determined up to addition of $a_8f_1$).
Altogether \eqref{c41} shows that
$$
d_4=6H_4^0(E,E,\bar{E},\bar{E})+3H_3^0(\tilde{r}_{200000}^0,\bar{E},\bar{E})+3H_3^0(\tilde{r}_{101000}^0-c_2f_2,E,\bar{E})+18H_3^0(E,\bar{E},f_2)^2,
$$
and the result of this calculation is given in equation \eqref{Formula for c4}.

\bibliographystyle{standard}

\end{document}